\documentclass[twocolumn,10pt]{asme2ej}

\usepackage{graphicx}
\usepackage{epsfig}
\usepackage{float}
\usepackage{subeqn}
\usepackage{mathrsfs}
\usepackage{indentfirst}
\usepackage{caption}
\usepackage{psfrag}
\usepackage{bm}
\usepackage{bbm}
\usepackage{color}
\usepackage{amsfonts}

\DeclareMathAlphabet{\mathcal}{OMS}{cmsy}{b}{n}  
\DeclareMathAlphabet{\mathcal}{OMS}{cmsy}{m}{n}

\newcommand{\Fig}[1]{Fig.~\ref{#1}}

\newcommand{\Ff}[1]{Figure~\ref{#1}}

\newcommand{\Eq}[1]{Eqn.~(\ref{#1})}

\newtheorem{definition}{Definition}

\newsavebox{\tablebox}

\title{Multi-objective Robust Optimization using a Post-optimality Sensitivity Analysis Technique: Application to a Wind Turbine Design}

\author{Weijun WANG$^{1}$, \,St\'{e}phane CARO$^2$\thanks{Corresponding author: Email: stephane.caro@irccyn.ec-nantes.fr, Tel: +33 (0)2 40 37 69 68, Fax: +33 (0)2 40 37 69 30}, \,Fouad BENNIS$^1$, \, Ricardo SOTO$^{3,4}$ , \, Broderick CRAWFORD$^{5,6}$
    \affiliation{$^1$  Ecole Centrale de Nantes,
    Institut de Recherche en Communications
    et Cybern\'etique de Nantes,\\
    1 rue de la No\"e, 44321
    Nantes, France\\
    Email: \{weijun.wang, fouad.bennis\}@irccyn.ec-nantes.fr \\

   $^2$ CNRS, Institut de Recherche en Communications
    et Cybern\'etique de Nantes,
    UMR CNRS 6597, France\\
    Email: stephane.caro@irccyn.ec-nantes.fr \\

	$^3$ Pontificia Universidad
Cat\'olica de Valpara{\'\i}so, \\Av. Brasil 2950, 2362807, Valpara{\'\i}so, Chile\\

 $^4$  Universidad Aut\'onoma de Chile, \\Av. Pedro de Valdivia 641, 7500138, Santiago, Chile\\

 $^5$  Universidad Finis Terrae, \\
Av. Pedro de Valdivia 1509, 7501015, Santiago, Chile\\

$^6$ Facultad de Ingenier\'ia y Tecnolog\'ia, Universidad San Sebasti\'an,\\ Bellavista 7, Recoleta, 8420524, Santiago, Chile\\

    Email: \{ricardo.soto, broderick.crawford\}@ucv.cl

    }
}

\begin{document}

\maketitle

\begin{abstract}
{\it
	Toward a multi-objective optimization robust problem, the variations in design variables and design environment parameters include the small variations and the large variations. The former have small effect on the performance functions and/or the constraints, and the latter refer to the ones that have large effect on the performance functions and/or the constraints. The robustness of performance functions is discussed in this paper. A post-optimality sensitivity analysis technique for multi-objective robust optimization problems is discussed and two robustness indices are introduced. The first one considers the robustness of the performance functions to small variations in the design variables and the design environment parameters. The second robustness index characterizes the robustness of the performance functions to large variations in the design environment parameters. It is based on the ability of a solution to maintain a good Pareto ranking for different design environment parameters due to large variations. The robustness of the solutions is treated as vectors in the robustness function space, which is defined by the two proposed robustness indices. As a result, the designer can compare the robustness of all Pareto optimal solutions and make a decision. Finally, two illustrative examples are given to highlight the contributions of this paper. The first example is about a numerical problem, whereas the second problem deals with the multi-objective robust optimization design of a floating wind turbine.
	}
\end{abstract}

\begin{nomenclature}

\entry{$\bf x$} {design variable vector}

\entry{$\bf p$}  {vector of design environment parameters}

\entry{$g_k$} {the $k$th constraint}

\entry{$\bf f$} {performance function vector}

\entry{$\mathcal{F}$} {feasible set}

\entry{$\mathcal{P}$} {Pareto optimal set}

\entry{$\sigma$} {standard deviation}

\entry{$\mu$} {expected value}

\entry{ $I_{RS}$} {robustness index against small variations}

\entry{ $I_F$} {feasibility index of a solution}

\entry{$I_P$} {Pareto optimality index of a solution}

\entry{ $h(\bf p)$}  {probability density function of $\bf p$}

\entry{$I_{rank}$} {individual's ranking}

\entry{$N$} {number of discrete values of design environment parameters}

\entry{ $I_{RL}$} {robustness index against large variations}

\entry{$\mathcal{P}_R(\mathcal{P})$} {the most robust solutions amongst the Pareto optimal solutions}

\entry{$P$}{produced power of wind turbine rotor}

\entry{$F_a$}{the thrust force in the partial load region of a wind turbine rotor}

\entry{$\gamma_r$}{the root twist angle}

\entry{$\gamma_t$}{the tip twist angle}

\entry{$c_r$ }{the chord length at the root}

\entry{$c_t$ }{the chord length at the tip}

\entry{$\omega$}{rotor rotational speed}

\entry{$r_{r}$}{root radius of the wind turbine}

\entry{$r_{t}$}{tip radius of the wind turbine}

\entry{$b$ }{number of blades}

\entry{$\rho$}{air density}

\entry{$v_{re}$} {reference wind speed}

\entry{HAWT} {Horizontal Axis Wind Turbine}

\entry{MOOP} {Multi-Objective Optimization Problem}

\entry{MOROP} {Multi-Objective Robust Optimization Problem}

\entry{RI}{Robustness Index}

\entry{RDP}{Robust Design Problem}

\entry{RF-Space}{Robustness Function Space}

\entry{DV} {Design Variable}

\entry{DEP} {Design Environment Parameter}

\entry{PF-Space}{Performance Function Space}

\entry{PDF}{Probability Density Function}

\entry{BEMT}{Blade Element Momentum Theory}

\entry{UD}{Uniform Distribution}

\entry{ND}{Normal Distribution}

\entry{N/A}{Not Affected}

\end{nomenclature}

\section{Introduction}  \label{Introduction}

Many design optimization problems are Multi-Objective Optimization Problems (MOOP) and are subject to uncertainties or variations in their parameters. Robustness is a product's ability to maintain its performance under the variations in its parameters. Robust design aims at maximizing product's robustness. In other words, it aims at minimizing the sensitivity of performance to variations without controlling the source of these variations~\cite{Caro2005}.
Sometimes, the robustness of a product is as important as or even more
important than its performance. So, to focus on the trade-off between robustness and performance of a product is meaningful. For such problems, namely Multi-Objective Robust Optimization Problems (MOROP), it is important to obtain design solutions that are both optimal and robust.

There exist some Robustness Indices~(RI) for MOROP in the literature. In fact, how to account for the variations in design parameters and how to measure robustness is a widely discussed problem~\cite{Beyer2007}.
Without loss of generality, the design parameters can be divided into two types, Design Variables~(DVs), which can be controlled by the designer, and Design Environment Parameters~(DEPs), which are uncontrollable parameters.

Existing works mainly focus on the small variations in DVs and DEPs~\cite{Deb2005,Deb2006,Saha2011,Gaspar-Cunha2006,Gaspar-Cunha2008,Barrico2006,Barrico2006a,Giassi2004,Caro2005,Caro2005a,Wang2012,Augusto2012a,Gunawan2005,Li2005,Li2006,Hu2009,Hung2011,Hung2013}. Moreover, some previous work focused on the variations in the performance
function model~\cite{Apley2006,Hassan2008,Hassan2009,Lu2009}. Toward MOOP, the designer's preference is also a type of variation in MOROP~\cite{Moshaiov2006,Avigad2011,Bui2012}.

Besides the small variations in DVs and DEPs, there may be large variations in DEPs in some engineering problems. For instance, the actual wind speed may vary greatly in a short time for a wind turbine rotor. In some areas, the environment temperature can fluctuate greatly between daytime and night. As a matter of fact, the problem with regard to large variations in DEPs can be considered as a special dynamic optimization problem, where the performance function and/or the constraints change with time~\cite{Cruz2011,Farina2003,Deb2007,Chen2011,Hassan2008,Hassan2009}.

One way to distinguish between small variations and large variations in design parameters is whether they are uncertain or not~\cite{Nazari2013}. Another way to distinguish between small and large variations is whether a linearization is reasonable over the range of variations. The distinction between small and large variations is commonly presented as local versus global sensitivity analysis in the literature on sensitivity analysis~\cite{Saltelli2008,Tomlin2011,Tong2008,Augustin2011}. The concept of large variations has been discussed in~\cite{Gunawan2004,Li2006,Lu2009,Lu2010} for comparing the difference between the linear and nonlinear performance functions. Nevertheless, in those discussions with regard to large variations, most of the existing criteria for MOROP aim at finding a design solution that gives desirable means, variances, quantiles or probabilities of violating constraints based on the distributions of the performance and constraint values with respect to the variations in the design parameters. Few criteria aim at finding a solution that maintains a good Pareto ranking for as many different DEP values in the sample space as possible. As a matter of fact, the Pareto ranking of a solution may be dramatically affected due to large variations in DEPs. In this paper, the distinction between small are large variations is made based on their effect on the performance functions. The small variations in design parameters refer to the ones that have small effect on the performance functions. The large variations in design parameters refer to the ones that have large effect on the performance functions.

Accordingly, a post-optimality sensitivity analysis technique for multi-objective robust optimization problems is discussed while considering both small and large variations in design parameters. Two robustness indices~(RI) are also introduced. The first one characterizes the robustness of MOOPs against small variations in the design parameters. It is based on the distributions of the performance function values with respect to these small variations. The second robustness index characterizes the robustness of MOOPs against large variations in DEPs. It is based on the solution's ability to maintain a good Pareto ranking for as many different design environments as possible. As a result, the designer can make a decision based on the robustness of the solutions.

The paper is organized as follows. Section~2 provides the theoretical background on Multi-Objective Optimization Problems~(MOOPs), Robust Design Problems~(RDPs) and Multi-Objective Robust Optimization Problems~(MOROPs). Two Robustness Indices~(RI) and the concept of Robustness Function Space (RF-Space) are introduced in Section~3.  Two illustrative examples are given in Section~4 to highlight the contributions of the paper. The first example is about a numerical problem, whereas the second problem deals with the multi-objective robust optimization design of a floating wind turbine. Conclusions and future work are presented in Section~5.

\section{Theoretical Background}
\label{sec:ProbFormul}
\subsection{Multi-Objective Optimization Problem}
 A general formulation of MOOP is given in \Eq {MOOP1}:
\begin{equation} \label{MOOP1}
    \begin{array}{rcl}
\textrm{minimize} & {\bf f}({\bf x}, {\bf p})=[f_1~ f_2 ~\ldots ~f_m]^T \\
\textrm{subject to}  & g_k({\bf x}, {\bf p})\leq 0, \quad k=1,\ldots,q \\
 & x_l^l \leq  x_l \leq  x_l^u, \quad l=1\ldots,n \\
    \end{array}
\end{equation}
where ${\bf f} ({\bf x}, {\bf p}) = [f_1~ f_2 ~\ldots ~f_m]^T $ denotes the $m$-dimensional vector of performance functions. ${\bf x} = [x_1~ x_2~ \ldots ~x_n]^T$ denotes the $n$-dimensional vector of DVs. Note that the nominal values of DV are controllable, $x_l^l$ and $x_l^u$ are the lower and upper bounds of $x_l$ respectively. ${\bf p} = [p_1 ~ p_2~ \ldots ~p_r]^T$ denotes the $r$-dimensional vector of DEP, which cannot be adjusted by the designer, and they are uncontrollable parameters. The functions $g_k({\bf x}, {\bf p})= [g_1~ g_2 ~\ldots ~g_q]^T$ are the constraints. A design that does not violate any constraint is called ``feasible''. The set of feasible solutions is called feasible set  and named $\mathcal{F}$.

Since there are some trade-offs amongst the $m$ conflicting objectives, none solution of $\mathcal{F}$ can dominate another solution of $\mathcal{F}$ . The optimization problem~(\ref{MOOP1}) generally has more than one optimal solutions. Those solutions are defined as Pareto optimal solutions, which cannot be dominated by any other feasible solution~\cite{Deb2001book,Li.X2009,Coello2006,Bui2012}. The set of all Pareto optimal solutions is called Pareto optimal set and named $\mathcal{P}$. The Pareto optimal solutions lie on a boundary in the Performance Function Space (PF-Space), called Pareto front.

\subsection{Robust Design Problem}

The concept of robustness is used in many fields such as engineering, biology, economy, computer science~\cite{Baker2008,Fink2009,Bui2012,Wieland2012,Kitano2004,Felix2008, Kuorikoski2007,Shahrokni2013}. In this paper, the robustness of a product is defined as follows:

\begin{definition} \label{def:Robustness0}
Robustness is a product's ability to maintain its performances under conditions of varying parameters.
\end{definition}

Robust design, a term originally introduced by Genichi Taguchi~\cite{Taguchi1993}, is a way of improving the quality of a product by minimizing the effect of variations, without eliminating the causes themselves. Many researchers refer to robustness~\cite{PARK2006,Arvidsson2008,Suh2001,Box1993,Caro2005,Zang2005a,Huang2006,Huang2006a,Lee2010}.

Although different expressions are used, their meanings are similar. In this paper, we define the robust design as follows:

\begin{definition} \label{def:Robust design}
  The design of a product is robust if and only if the performances of the good under design is as little sensitive as possible to variations and uncertainties.
\end{definition}

\subsection{Multi-Objective Robust Optimization Problem}

A MOROP aims to find out a solution that is feasible, optimal and robust. To come up with an feasible, optimal and robust design, the three following scenarios have been identified.

\begin{enumerate}
	\item \textbf{Optimal design and robust design are equally important:}
 It means that performance functions and robustness functions are optimized simultaneously. The designer can use the new performance functions by adding the effects of robustness, instead of original performance functions~\cite{Deb2005, Deb2006,Gaspar-Cunha2006,Gaspar-Cunha2008,Ferreira2008}. Moreover, the robustness functions can be considered as new constraints~\cite{Barrico2006,Barrico2006a,Barrico2007,Barrico2009,Gunawan2005,Li2006,Hu2009,Hu2011b} or additional optimization objectives~\cite{Li2005,Nejlaoui2011,Nejlaoui2012,Erfani2010,Erfani2012,AitBrik2006,Ghanmi2007,Ghanmi2011}.
	\item \textbf{Optimal design is primary and robust design is secondary:} It is one post-optimality approach. The final solution is selected from the Pareto optimal set based on the robustness criterion~\cite{Forouraghi2000,Olvander2005,Avigad2005,Moshaiov2006,Avigad2011,Xu2011,Hung2011,Hung2013,Cromvik2011}.
	\item \textbf{Robust design is primary and optimal design is secondary:} The final solution is selected from the most robust solutions based on the values of the performance functions~\cite{Wang2012}.
\end{enumerate}

In this paper, the optimal design is supposed to be of primary importance and the robust design of secondary importance. Although there exist some papers dealing with this method ~\cite{Forouraghi2000,Olvander2005,Avigad2005,Moshaiov2006,Avigad2011,Xu2011,Hung2011,Hung2013,Cromvik2011}, which is also called as post-optimality sensitivity analysis, the final solution is selected from the Pareto optimal set based on its robustness. Moreover, in this paper we consider both small variations in design parameters and large variations in DEPs simultaneously.

\section{New Robustness Indices for MOROP} \label{Proposed_method}

Two Robustness Indices~(RI) are introduced in this section. Then, the concept of RF-Space is presented and the solution robustness is defined as a vector in the RF-Space. Finally, a discussion on this new method is presented.

\subsection{Robustness Index with regard to the Small Variations in DVs and DEPs}

\begin{figure}[!htbp]
\begin{center}
\includegraphics[width=0.48\textwidth]{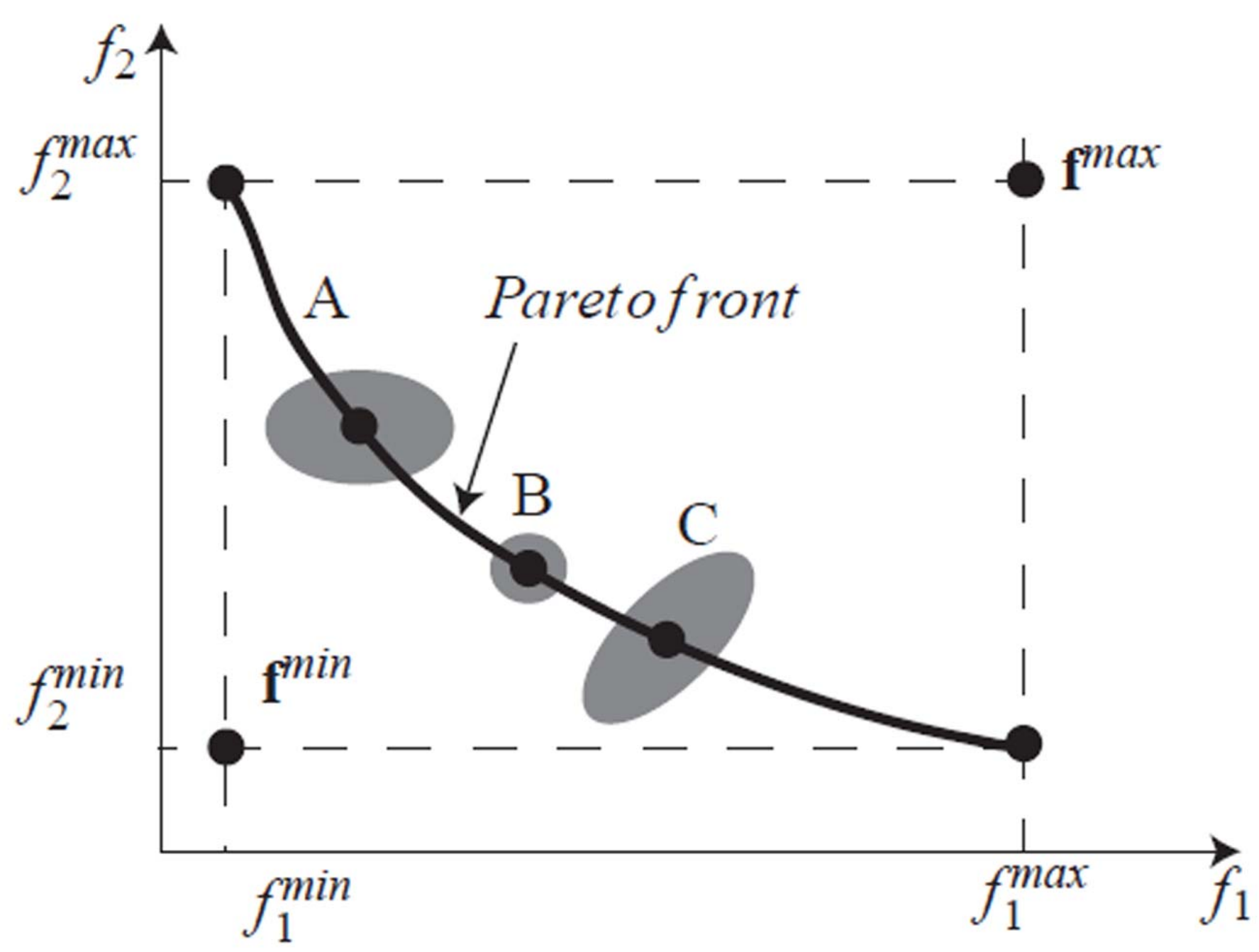}
\caption{The performances of solutions vary around their nominal values in the PF-Space, with regard to small variations in DVs and DEPs}
\label{IR_S_MOOP}
\end{center}
\end{figure}

\Ff {IR_S_MOOP} illustrates the solutions of a simple multi-objective optimization problem. It has two performance functions ($f_1$,$f_2$,), some DVs ($\bf x$) and DEPs ($\bf p$). Let us select solution~A, solution~B and solution~C, which are located on the Pareto front. When there are small variations in DVs and DEPs, then their performances varies in the grey area around their nominal values in the PF-Space, as shown in~\Fig {IR_S_MOOP}. Obviously, in the PF-Space, solution~B has smaller variations around its nominal values than solutions~A and~C. It means than solution~B is more robust than solutions~A and~C. However, the size of the grey areas associated with solution A and solution C are the same, but their shapes are different. So the comparison between the robustness of solution A and solution C is not an easy task. Toward the different types of uncertainties, there are different methods to define RI for MOOP, such as using the worst case scenario~\cite{Gunawan2005,Li2005,Li2006,Li2008,Hu2009}, the expectancy measure ~\cite{Gaspar-Cunha2006,Gaspar-Cunha2008,Deb2005,Deb2006,Saha2011,Barrico2006,Barrico2006a,PARK2006,Beyer2007,Du2009}, and the probabilistic threshold method~\cite{Barrico2006,Barrico2006a,Beyer2007}.

\begin{figure}[!htbp]
\begin{center}
\includegraphics[width=0.48\textwidth]{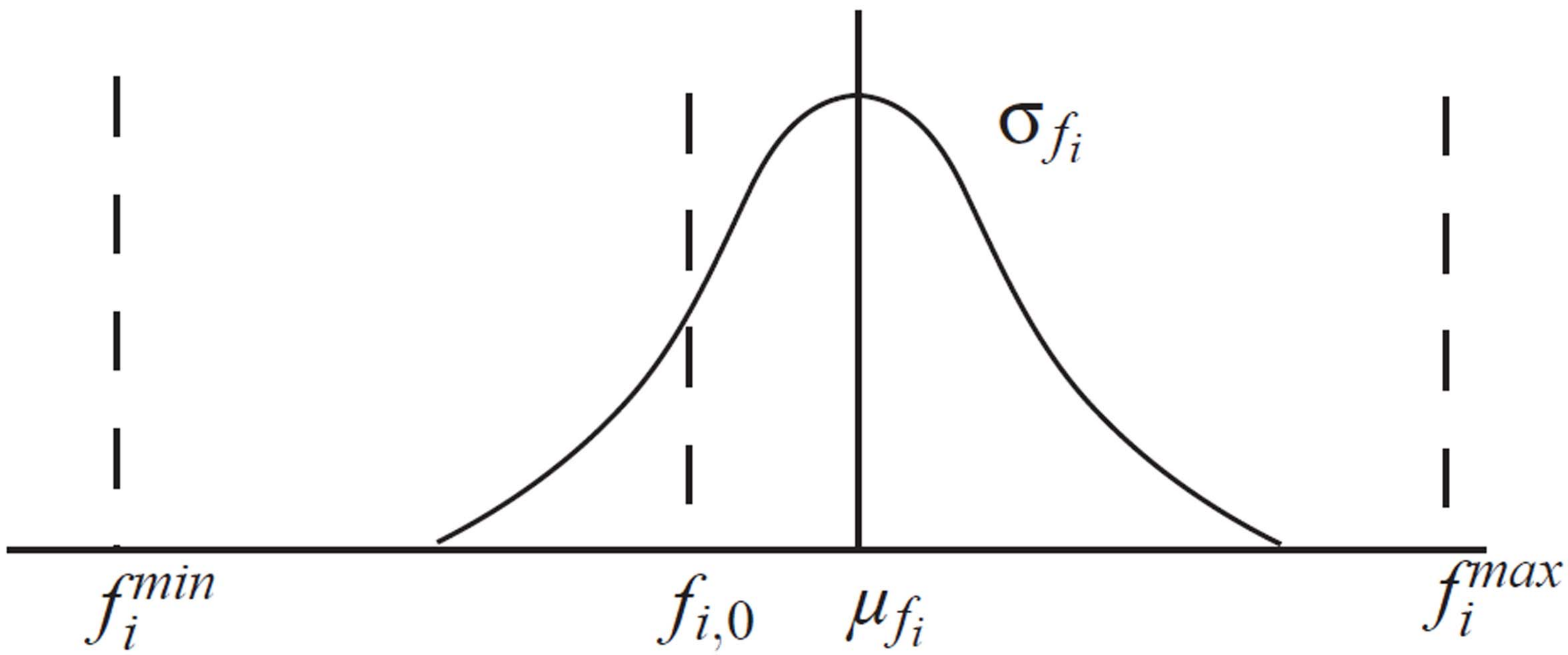}
\caption{The distribution of the $i$th performance function as a function of small variations in DVs and DEPs}
\label{IR_S_1}
\end{center}
\end{figure}

We select a RI against the small variations in~DVs and DEPs, based on the expectancy measure. A new RI ($I_{RS}$) is introduced and is based on the actual variations of performances. Under the small variations in DVs and DEPs, the actual performances are distributed around the nominal values. In the PF-Space, the actual performances follow a multivariate distribution. Here, in order to simplify the index associated with each performance function, the standard deviation~$\sigma$ of the actual performances is used as a measure for the robustness: the smaller the standard deviation, the more robust the design. The absolute value of the difference between the expected value~($\mu$) and the nominal value~($f_0$) is also a robustness measure: the smaller the absolute difference, the more robust the design. \Ff{IR_S_1} shows an example,  which follows a normal distribution, $\sigma_{f_i}$ and $\mu_{f_i}$ are the standard deviation and expected value of the $i$th performance function under small uncertainties; $f_{i,0}$ is the nominal value; $f^{max}_i$ and $f^{min}_i$ are the maximum and the minimum values of the $i$th performance function on the Pareto front, respectively.

Here, we assume that the standard deviation ($\sigma$) and the absolute value of the difference between~$\mu$ and~$f_0$ have the same importance for the designer. Moreover, we assume that the robustness of all performance functions has the same importance. The robustness index $I_{RS}$ of a MOOP with respect to small variations in DVs and DEPs is defined as a scalar. The robustness of each performance function is also normalized. To normalize the sum of the standard deviation ($\sigma_{f_i}$) and the absolute difference ($|\mu_{f_i}-f_{i,0}|$), we divide it by the difference between the two extreme values of the $i$th performance function, namely, $f^{max}_i-f^{min}_i$. As a result, $I_{RS}$ is defined as follows:
\begin{eqnarray} \label{IRS}
	I_{RS}({\bf x}) = \sqrt{ \sum\limits_{i=1}^{m} \left ( \frac{\sigma_{f_i}+|\mu_{f_i}-f_{i,0}|}{f^{max}_i-f^{min}_i} \right )^2}
\end{eqnarray}
The smaller~$I_{RS}$, the more robust the design.

Note that even if the variations are small, the constraints may not be satisfied due to variations. This refers to the robustness of constraints (reliability). As a matter of fact, the index~$I_{RS}$ is not discussed deeply in this paper. It is just based on the distributions of the performance function values with respect to these small variations. To simplify the index, the reliability with regard to small variations is not considered in this paper.

\subsection{Robustness Index with regard to the Large Variations in DEPs}

Since the DVs are controllable and the DEPs are uncontrollable, we assume that there are only large variations in DEPs. The mapping functions between the DVs and performance functions can change a lot due to large variations in DEPs.
For a MOOP while considering large variations in DEPs, we assume that the initial DEPs are ${\bf p}={\bf p}_0$, the set of feasible solutions is named $\mathcal{F}_{0}$ and the set of Pareto optimal solutions is denoted $\mathcal{P}_{0}$. It is noteworthy that those Pareto optimal solutions are alternative solutions for the designer.

Since large variations in DEPs exist, the DEPs may change from ${\bf p}_0$ to ${\bf p}_{new}$. The design environment parameters~$\bf p$ are supposed to take $N$~discrete values: ${\bf p}_1, {\bf p}_2, \cdots, {\bf p}_N$. The Probability Density Functions (PDF) of $\bf p$ are $h({\bf p}_1),h({\bf p}_2), \cdots, h({\bf p}_N)$. The initial DEPs ${\bf p}_0$ are equal to the ones having the maximum PDF amongst the $N$ discrete values.
The corresponding feasible sets are $\mathcal{F}_{1}, \mathcal{F}_{2}, \cdots, \mathcal{F}_{N}$. A feasible and Pareto optimal solution in $\mathcal{P}_0$ may not be Pareto optimal, and not feasible in the new environments. To compare the solution's robustness against large variations in DEPs, the traditional methods are not applicable. Toward a MOOP, a definition of the solution's robustness against large variations in DEPs is proposed thereafter:
\begin{definition} \label{def:Robustness}
  Toward a multi-objective optimization problem against large variations in DEPs, solution's robustness against large variations in DEPs is a measure of its ability to be optimal in different design environments.
\end{definition}

\begin{figure}[!htbp]
\begin{center}
\includegraphics[width=0.48\textwidth]{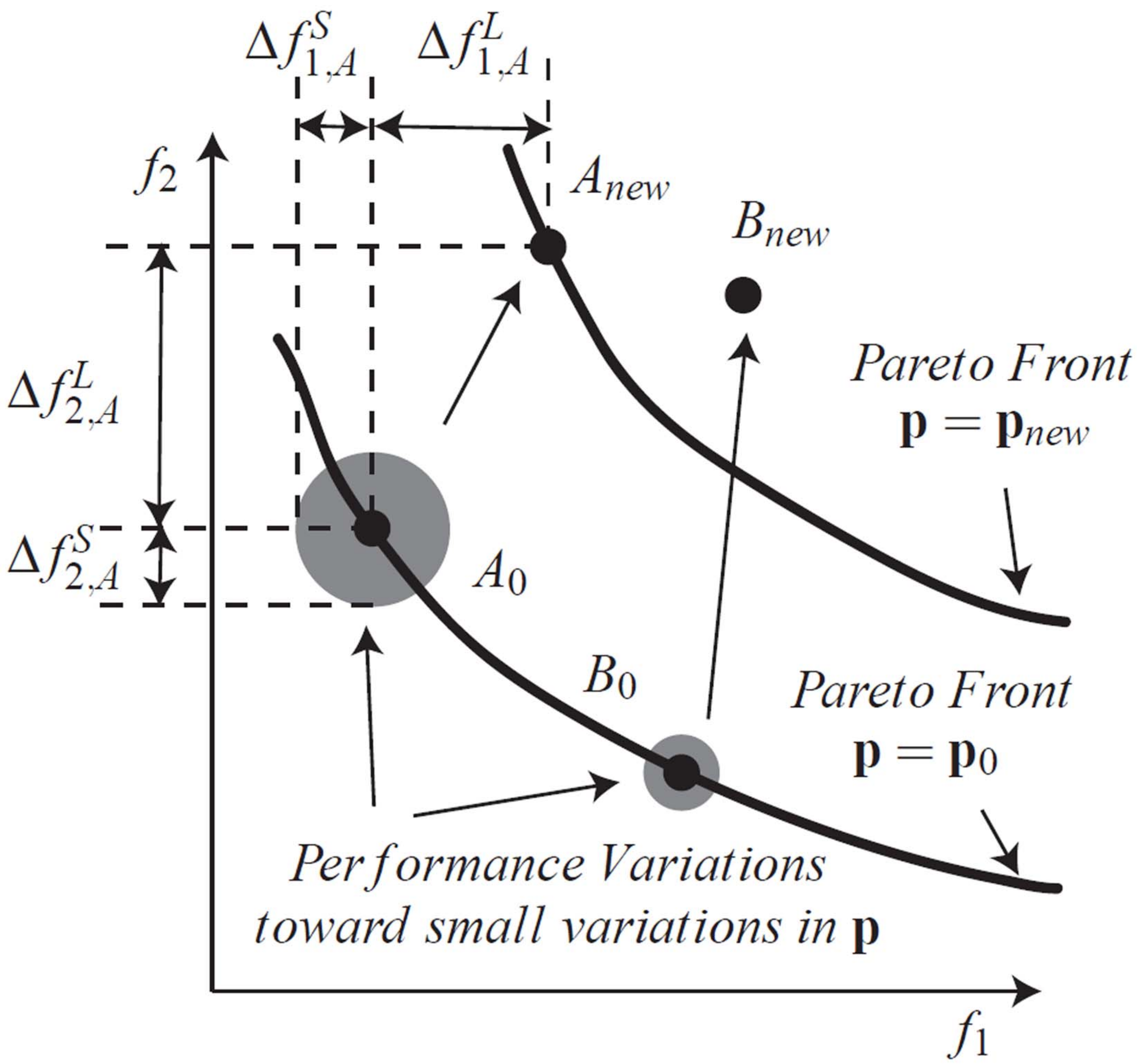}
\caption{The performances of solutions vary greatly in the PF-Space due to large variations in DEPs}
\label{RL-DEP-change-1}
\end{center}
\end{figure}

\Ff {RL-DEP-change-1} illustrates the solutions of a simple MOOP. Similarly to \Fig {IR_S_MOOP}, the problem has two objective functions ($f_1$,$f_2$), some DVs ($\bf x$) and DEPs ($\bf p$). ${\bf p}_0$ is the nominal values of DEPs. Let us consider Pareto-optimal solutions~$A$ and~$B$ as examples. The grey area shows the variations in performances of solution $A$ and solution $B$ in the PF-Space, under small variations in DVs and DEPs. As seen before, we can conclude that solution $B$ is more robust than solution $A$ with regard to small variations in DVs and DEPs.

On the contrary, it may not be the case if there are large variations in DEPs. As shown in \Fig {RL-DEP-change-1}, in a new environment, we assume that when the DEPs ${\bf p}$ change from ${\bf p}_0$ to ${\bf p}_{new}$, there are large variations in the PF-Space for all the solutions. For instance, in the PF-Space, solution A moves from $A_0$ to $A_{new}$, solution B moves from $B_0$ to $B_{new}$. Both solutions have quite large variations in $\bf f$. Then, the following question remains: \emph{``How can we compare the robustness of solutions~$A$ and~$B$"?}.

As shown in \Fig {RL-DEP-change-1}, toward solution $A$, $\Delta f_{1,A}^S$ and $\Delta f_{2,A}^S$ represent the largest distance of the actual performance function values and nominal performance function values in $f_1$ and $f_2$ respectively, with regard to the small variations in DVs and DEPs. $\Delta f_{1,A}^L$ and $\Delta f_{2,A}^L$ represent the largest distance of the actual performance function values (assuming only two samples: $A_0$ and $A_{new}$) and nominal performance function values in $f_1$ and $f_2$ respectively, with regard to large variations in DEPs.
If $\Delta f_{1,A}^L>> \Delta f_{1,A}^S$, $\Delta f_{2,A}^L>> \Delta f_{2,A}^S$ , and $\Delta f_{1,B}^L>> \Delta f_{1,B}^S$, $\Delta f_{2,B}^L>> \Delta f_{2,B}^S$, we can see that the caused actual performance function values of both solution $A$ and solution $B$ are sufficiently far from their nominal ones, with regard to large variations in DEPs. As a consequence, we can conclude that the traditional methods, which provide some results based on the difference between actual performance function values and nominal ones, make little sense. Therefore, the method proposed in this paper aims to compare their relative positions in the PF-Space associated with the new environment. From \Ff {RL-DEP-change-1}, we can see that, in the new environment, ${\bf p}={\bf p}_{new}$, solution $A$ is still on the new Pareto front, which means that it is still one of the best solutions for the designer. Meanwhile, solution $B$ is far away from the new Pareto front, which means that it is no longer a good choice for the designer. As a result, we can conclude that solution $A$ is more robust than solution $B$ with regard to the large variations in DEPs.

\begin{figure}[!htbp]
\begin{center}
\includegraphics[width=0.48\textwidth]{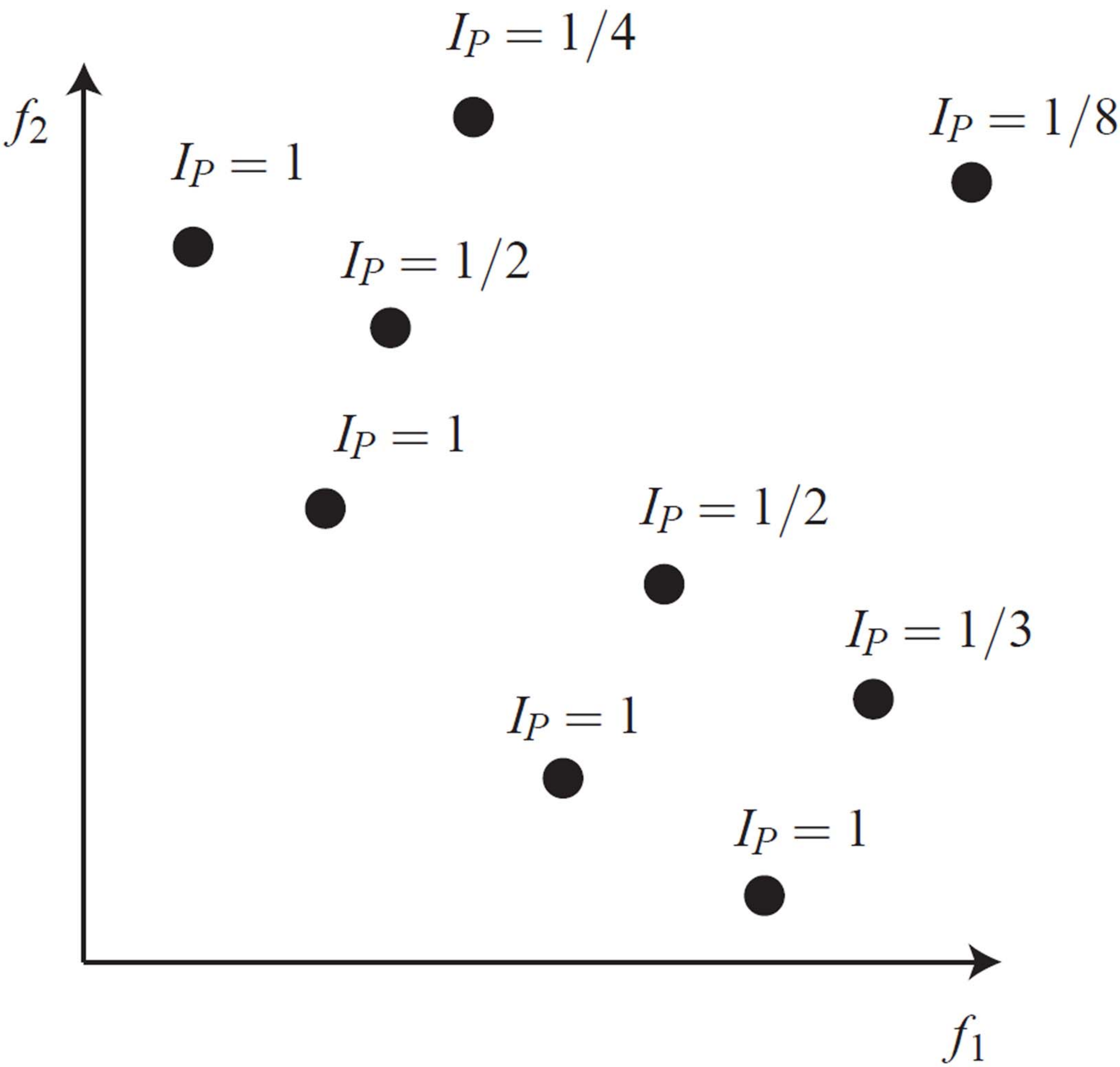}
\caption{The positions of the Pareto optimal solutions in a new environment}
\label{RL-Rank}
\end{center}
\end{figure}

As a result, a RI with regard to large variations in DEPs is defined thereafter. Toward the discrete probability distribution of DEPs, mathematically, the RI $I_{RL}$ of a solution~${\bf x}$ is defined as follows:
\begin{subequations} \label{Robustness_index}
\begin{eqnarray}
I_{RL}({\bf x}) = 1-I_F({\bf x}) \sum\limits_{1}^{N} I_P({\bf x},{\bf p}_j)h({\bf p}_j), \label{Robustness_indexIRL} \\
 I_F({\bf x}) = \left\{
  \begin{array}{ll}
    1, &  \hbox{$ \forall j=1,2,\cdots ,N, {\bf x} \in \mathcal{F}_j $;} \\
    0, &  \hbox{$ \exists j=1,2,\cdots ,N, {\bf x} \not\in \mathcal{F}_j$.}
  \end{array}
\right. \label{Robustness_indexIF}
\\
I_P({\bf x},{\bf p}_j)=1/I_{rank}({\bf x},{\bf p}_j) \label{Robustness_indexIP}
\end{eqnarray}
\end{subequations}
where $I_F$ is the Feasibility Index of the solution; $I_P$ is the Pareto optimality Index of the solution; $h(\bf p)$ is the PDF of $\bf p$; $I_{rank}({\bf x},{\bf p}_j)$ is the individual's ranking in the new environment where ${\bf p}={\bf p}_{j}$ and amounts to the number of individuals by which it is dominated amongst the alternative solutions, plus one~\cite{Fonseca1998,Augusto2006}. For a better understanding, a simple example is shown in \Fig {RL-Rank}. In a new environment, if a solution is still non dominated by any other solution, then the $I_P$ value will be equal to one for that solution. Otherwise, the $I_P$ value will be lower than one, but greater than zero. Note that the number of the alternative solutions affects the value. However, the proposed definition can divide the alternative solutions into different groups based on their robustness.

Note that the index $I_{RL}$ is bounded between zero and one. On the one hand, if a solution is feasible in all environments and cannot be dominated by any other solution in all possible environments, then $I_{RL}=0$. On the other hand, if a solution is non-feasible in some new environments, then $I_{RL}=1$. If each solution belongs to the set of Pareto optimal solutions $\mathcal{P}_0$, then if it is feasible in all environments and its $I_{RL}$ value will be greater than or equal to zero and smaller than  one.

Note that the individual's ranking~$I_{rank}({\bf x},{\bf p}_j)$ corresponds to a specific value of DEPs: ${\bf p}_j$, where~$j=1,2,\cdots ,N$. In case there exist continuous probability distributions of the DEPs, $N$ becomes infinite and it is difficult to assess the robustness index $I_{RL}$ for a solution~${\bf x}$. However, since the domain of the DEPs can be partitioned into many small parts, we can simplify such a problem by using a discrete probability distribution of the DEPs.

Finally, the smaller $I_{RL}$, the more robust the design with regard to large variations in  DEPs.

\subsection{Robustness Function Space}

\begin{figure}[t]
\begin{center}
\includegraphics[width=0.48\textwidth]{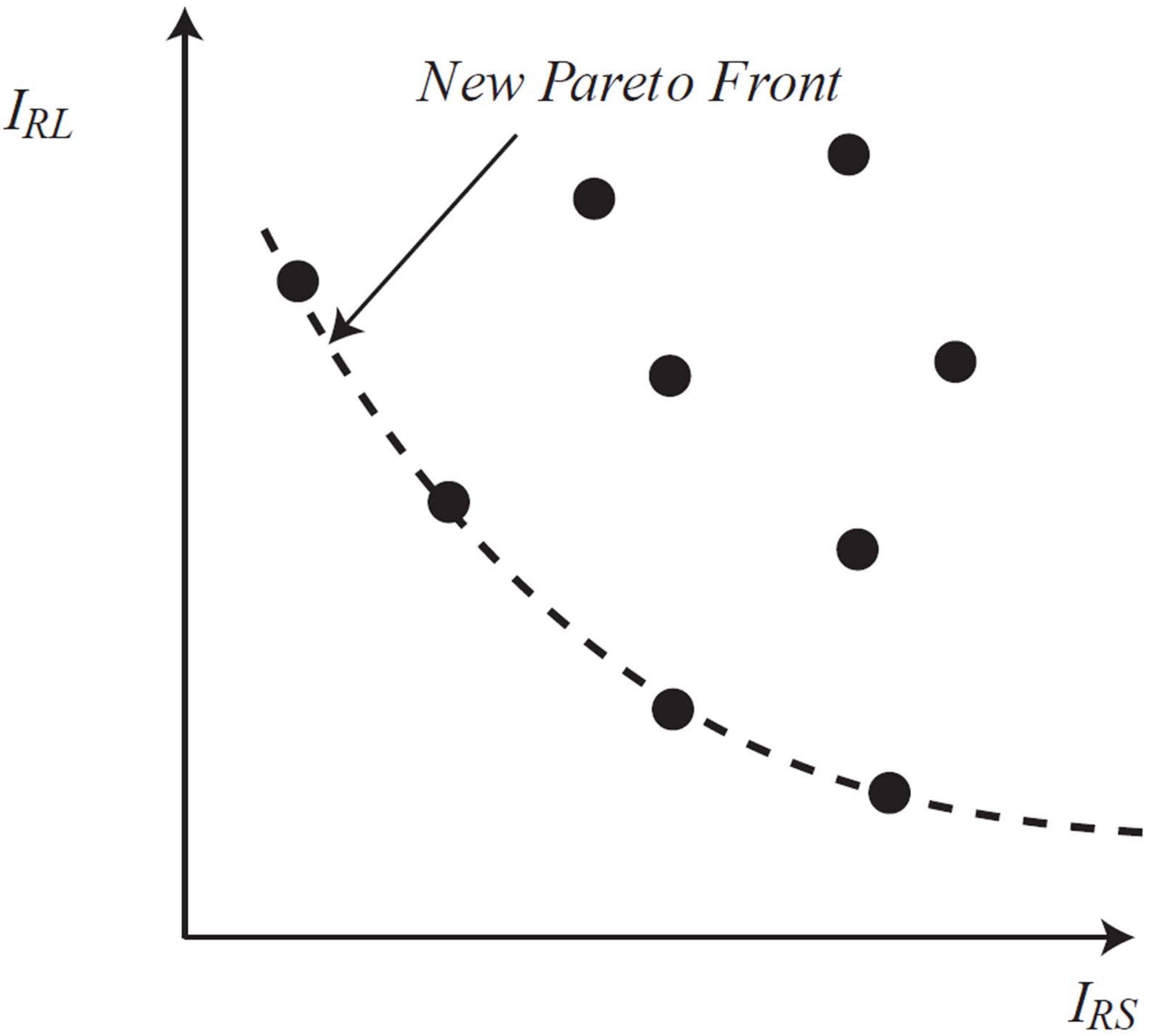}
\caption{Each Pareto optimal solution has a corresponding position in the RF-Space}
\label{RFS-1}
\end{center}
\end{figure}

In this paper, we consider not only the RI against small variations, $I_{RS}$, but also the RI against large variations, $I_{RL}$. The RI of a Pareto-optimal solution is represented as a vector. The designer can analyze the robustness of Pareto optimal solutions in the Robustness Function Space~(RF-Space). For a better understanding, we take $I_{RS}$ as one dimension of the RF-Space, another dimension of the RF-Space is $I_{RL}$, as shown in~\Fig {RFS-1}.

Thanks to the proposed method, each Pareto optimal solution has a corresponding position in the RF-Space, as shown in \Fig {RFS-1}. If $I_{RS}$ and $I_{RL}$ are not conflicting, then the designer will be able to select the most robust solution immediately, namely, the solution that minimizes both $I_{RS}$ and $I_{RL}$. If $I_{RS}$ and $I_{RL}$ are two conflicting objectives, then a new Pareto front in the RF-Space will appear. The Pareto-robust solutions represent the most robust solutions amongst the Pareto optimal solutions. The set of Pareto-robust solutions is named $\mathcal{P}_R(\mathcal{P})$. Finally, the designer can select the final solution from this set according to his/her requirements.

\subsection{Flow chart}

\begin{figure}[!htbp]
\begin{center}
\includegraphics[width=0.48\textwidth]{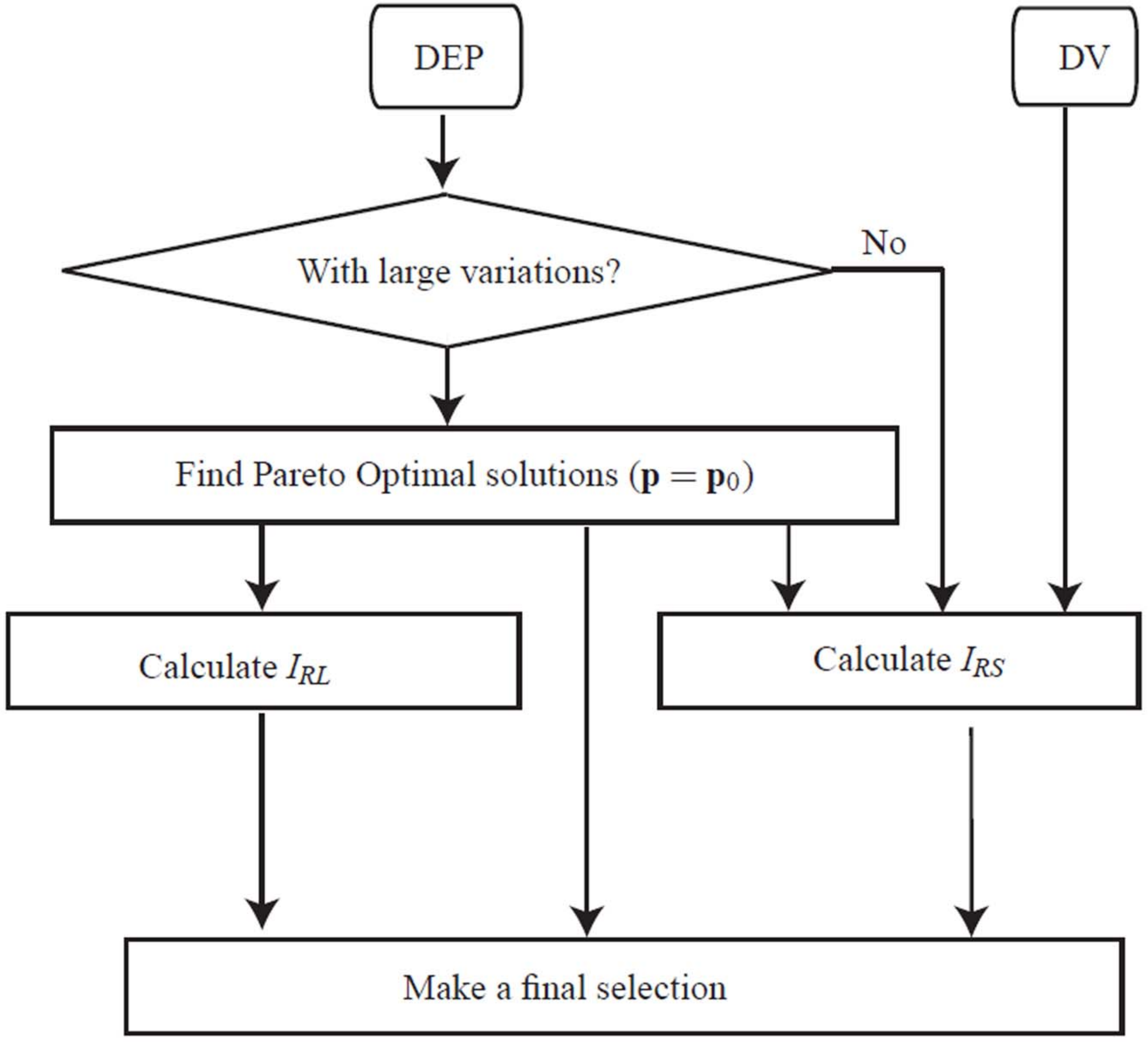}
\caption{A flow chart illustrating the proposed post-optimality sensitivity analysis technique}
\label{Flow-chart1}
\end{center}
\end{figure}

A flow chart is illustrated in \Ff {Flow-chart1} for a better understanding of the proposed post-optimality sensitivity analysis technique. The first step is to determine whether the design environment parameters~(DEP) are subject to large variations or not.  The distinction between small and large variations is made with regard to their effect on the performance functions. The small variations in design parameters refer to the ones that have a small effect on the performance functions. The large variations in design parameters refer to the ones that have large effect on the performance functions. In the second step, Pareto optimal solutions are obtained for different values of DEP. Those different values of DEP are due to large variations in DEP. In the third step, robustness indices $I_{RS}$ and $I_{RL}$ are calculated for each solution. The new Pareto ranking of the obtained Pareto optimal solutions are obtained in different design environments. Finally, the designer can make a decision from the obtained Pareto optimal solutions based on their robustness, namely, based on the values of $I_{RS}$ and $I_{RL}$ for each solution.

\section{Examples}

In this section, two illustrative examples are given in order to highlight the contributions of the proposed approach. The first example is a simple numerical example and the second example is a MOROP of a wind turbine blade design, in which the variations include not only the small variations in DVs and DEPs, but also large variations in DEPs.

\subsection{Numerical illustrative example}
The problem is defined in \Eq {Numerical_eg1}. There are two performance functions $f_1$ and $f_2$, one design variable~$x$, one design environment parameter~$p$ and one constraint:
\begin{subequations} \label{Numerical_eg1}
\begin{eqnarray}
  \textrm{minimize } &   f_1(x,p)=x+p/2\\
  \textrm{minimize } &   f_2(x,p)=(x-p)^2\\
  \textrm{subject to}  & 1 \leq x \leq 10 \\
  ~                    & f_1(x,p) \geq 3
\end{eqnarray}
\end{subequations}
$x$~is subject to small variations and the variations in $x$ follow a uniform distribution and are bounded between~$-0.1$ and $0.1$. $p$~is subject to large variations and can take three possible values: $p = 3$,  $p = 5$ and $p = 8$. The probabilities for $p$ to take those three values are $h(1) =0.2$, $h(2) = 0.5$ and $h(3) = 0.3$, respectively. The initial value of $p$ is equal to~5.

\begin{figure}[!htbp]
\begin{center}
\scriptsize
\includegraphics[width=0.48\textwidth]{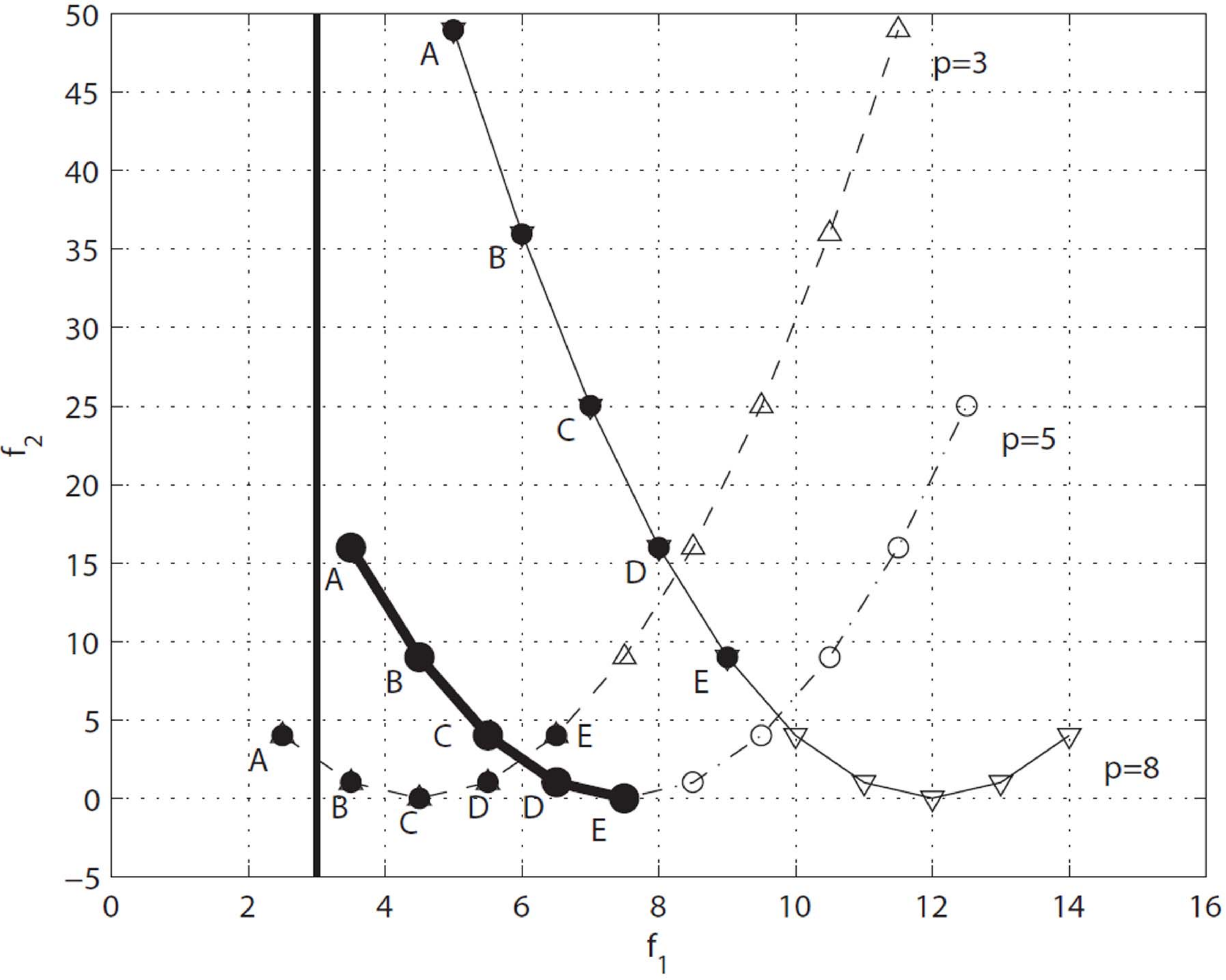}
\caption{Values of the performance functions~$f_1$ and~$f_2$ for $p = 3$, $p = 5$ and $p = 8$, respectively}
\label{numerial_example_1}
\end{center}
\end{figure}

\Ff {numerial_example_1} shows the values of the performance functions~$f_1$ and~$f_2$ for $p = 3$, $p = 5$ and $p = 8$, respectively, i.e., in different design environments. For $p=5$, let us consider the five Pareto optimal solutions: solution A obtained for $x=1$; solution B obtained for $x=2$; solution C obtained for $x=3$; solution D obtained for $x=4$ and solution E obtained for $x=5$. These solutions form the Pareto front in the PF-Space and are alternative solutions for the designer. Their new positions for $p=3$ and $p=8$ in the PF-Space are shown in \Fig{numerial_example_1}. The vertical line represents the constraint related to the first performance function. Note that the solutions to the right of this line are feasible, whereas the solutions to the left are not feasible.

To assess the robustness index $I_{RS}$ for each alternative solution, the samples around the nominal values are generated by Latin Hypercube Sampling (LHS)~\cite{McKay1979,Luo2008a,Petrone2011}. 1000 samples for each alternative solution are generated with regard to the small variations in DV. Then $I_{RS}$ for each alternative solution can be assessed with~\Eq {IRS}.

To determine the robustness index $I_{RL}$ associated with each alternative solution, their new positions should be located when $p$ changes from its initial value to its new value. The results can be found from \Fig {numerial_example_1}. When $p=3$, solution A becomes unfeasible, then the index $I_F$ defined in \Eq {Robustness_indexIF} is null for solution A. Solution~B and solution~C become Pareto optimal amongst these five solutions. Solution~D and Solution~E are dominated by other alternative solutions. When $p=8$, all the alternative solutions are feasible and can not be dominated by others. Then $I_{RL}$ for each alternative solution can be calculated with \Eq{Robustness_index}.

\begin{table}[!htbp]
\caption{Comparison of the alternative solutions with regard to variations in small variations in DV~$x$ and large variations in DEP~$p$}
\label{table5solutions_Numericaleg}
\begin{center}
\small
\renewcommand{\arraystretch}{2}
\begin{tabular}{ c c c c c c c}
 \hline
  \hline
  Solution & $f_1$& $f_2$ & $\Delta f_1^S$ &  $\Delta f_2^S $ & $\Delta f_1^L$& $\Delta f_2^L$ \\
  \hline
 A  &3.5 &	16 & 0.05&	0.4024  &1.5   &	33  \\
 B  &4.5 &	9 &	0.05 &	0.3024  &1.5   &27  \\
 C  &5.5 &	4 & 0.05 & 0.2024  &	1.5   &21  \\
 D  &6.5 &	1 &	0.05 &	0.1025  &1.5   &15  \\
 E  &7.5 &	0 &	0.05 &	0.0025 &	1.5   &9   \\
  \hline
  \hline
\end{tabular}
\end{center}
\end{table}

Note that all Pareto optimal solutions have sufficiently large variations in performance functions with regard to large variations in DEPs. In Tab.~\ref{table5solutions_Numericaleg}, $f_1$ and $f_2$ denote the nominal values of the performance functions for the alternative solutions. $\Delta f_1^S$ and $\Delta f_2^S$ represent the largest distance of the actual performance function values (1000 samples for each alternative solution) and nominal performance function values in $f_1$ and $f_2$ respectively, with regard to the small variations in DV. $\Delta f_1^L$ and $\Delta f_2^L$ represent the largest distance of the actual performance function values (3 samples for each alternative solution) and nominal performance function values in $f_1$ and $f_2$ respectively, with regard to large variations in DEP.
If the traditional methods (results based on the difference between actual performance function values and nominal ones) are selected to measure the robustness of the alternative solutions, here is the order of the six solutions from the most robust one to the least robust one: (1)~E; (2)~D; (3)~C; (4)~B; (5)~A. Even if solution~E is the most robust one, $\Delta f_{1,E}^L$ and $\Delta f_{2,E}^L$ are still quite large. Obviously, $\Delta f_{1,E}^L>> \Delta f_{1,E}^S$ and $\Delta f_{2,E}^L>> \Delta f_{2,E}^S$, so with regard to large variations in DEPs, the traditional methods make little sense.

Here is the order of the six solutions from the most robust one to the least robust one while using the robustness index~$I_{RL}$ defined in Eq.~(\ref{Robustness_indexIRL}): (1)~B; (1)~C; (3)~D; (4)~E; (5)~A. Solution B and solution C are Pareto optimal  no matter the value of DEP~$p$. Therefore, the designer may select solution~B or solution~C instead of solution~E.

\begin{figure}[!htbp]
\begin{center}
\scriptsize
\includegraphics[width=0.48\textwidth]{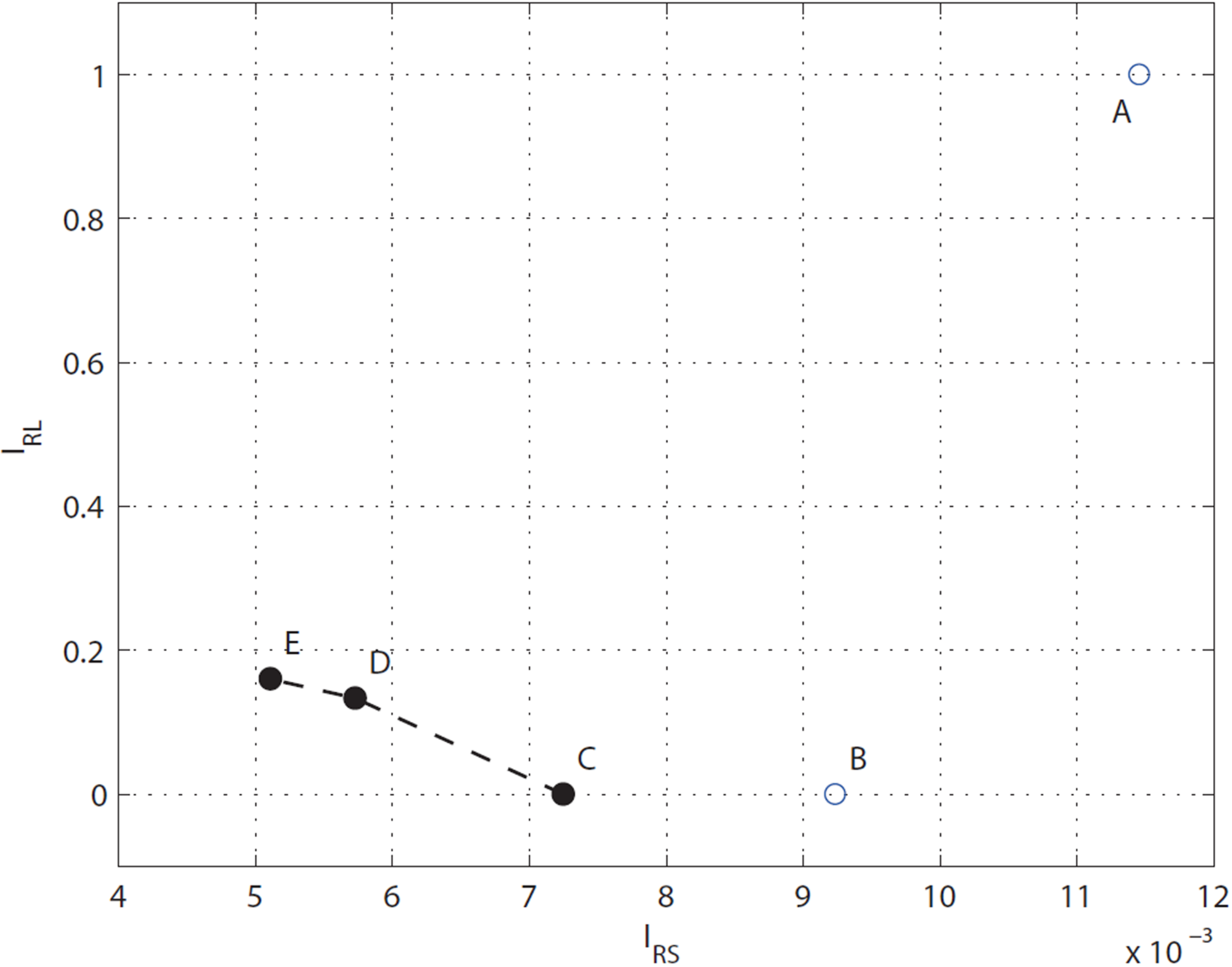}
\caption{The positions of the alternative solutions for the numerical example in the RF-Space}
\label{numerial_example_2}
\end{center}
\end{figure}

\Ff {numerial_example_2} illustrates the positions of the alternative solutions in the RF-Space. In this space, one dimension is $I_{RS}$, another dimension is $I_{RL}$. In this example, $I_{RS}$ and $I_{RL}$ are two conflicting objectives, as shown in \Fig {numerial_example_2}. Solutions C, D, E form a new Pareto front in the RF-Space, which represents the most robust solutions amongst these five alternative solutions.  For a better understanding, the performance function values and their RIs are given in Tab.~\ref {table5solutions_Numericaleg}. The indices $I_F$ and $I_P$ are also given in Tab. \ref {table5solutions_Numericaleg}.

\begin{table*}[!htbp]
\caption{Comparison of five alternative solutions for the numerical example}
\label{table5solutions_Numericaleg}
\begin{center}
\small
\renewcommand{\arraystretch}{2}
\begin{tabular}{ c c c c c c c c c c }
 \hline
  \hline
  Solutions & $f_1$ & $f_2$ & $I_{RS}$ & $I_F(x)$ & $I_P(x,p_1)$ & $I_P(x,p_2) $ & $I_P(x,p_3)$ & $I_{RL}$ \\
  \hline
 A  &3.5 &	16 & 0.0115&	0  &	1          &	1   &	1 &	1 \\
 B  &4.5 &	9 &	0.0092 &	1  &	1           &1   &	1 &	0 \\
 C  &5.5 &	4 &	0.0072 & 1  &	1           &1   &	1 &	0 \\
 D  &6.5 &	1 &	0.0057 &	1  &	0.3333  &1   &	1 &	0.1333 \\
 E  &7.5 &	0 &	0.0051 &	1  &	0.2        &1   &	1 &	0.16 \\
  \hline
  \hline
\end{tabular}
\end{center}
\end{table*}

Thanks to the proposed method, the designer can select the final solution from the new Pareto front in the RF-Space, according to his/her requirement. For example, if the designer prefers the $I_{RS}$ index, solution E will be selected. If the designer prefers the $I_{RL}$ index, solution C will be selected.


\subsection{Multi-Objective Robust Optimization Design of a Floating Wind Turbine}

\subsubsection{Problem Formulation}

\Ff {wind_turbine} illustrates a schematic of a floating horizontal wind turbine rotor with two simplified morphing blades. Each blade can adjust its tip twist angle and root twist angle according to the reference wind speed~\cite{Wang2012a}. However, it is difficult to adjust the twist angles at all time. Therefore, the twist angle is assumed to be adjusted according to the average wind speed. After setting the tip twist angle ($\gamma_t$) and root twist angle ($\gamma_t$), the twist angles of the other elements are adjusted automatically and they are linearly distributed along the blade.

\begin{figure}[!htbp]
\begin{center}
\includegraphics[width=0.48\textwidth]{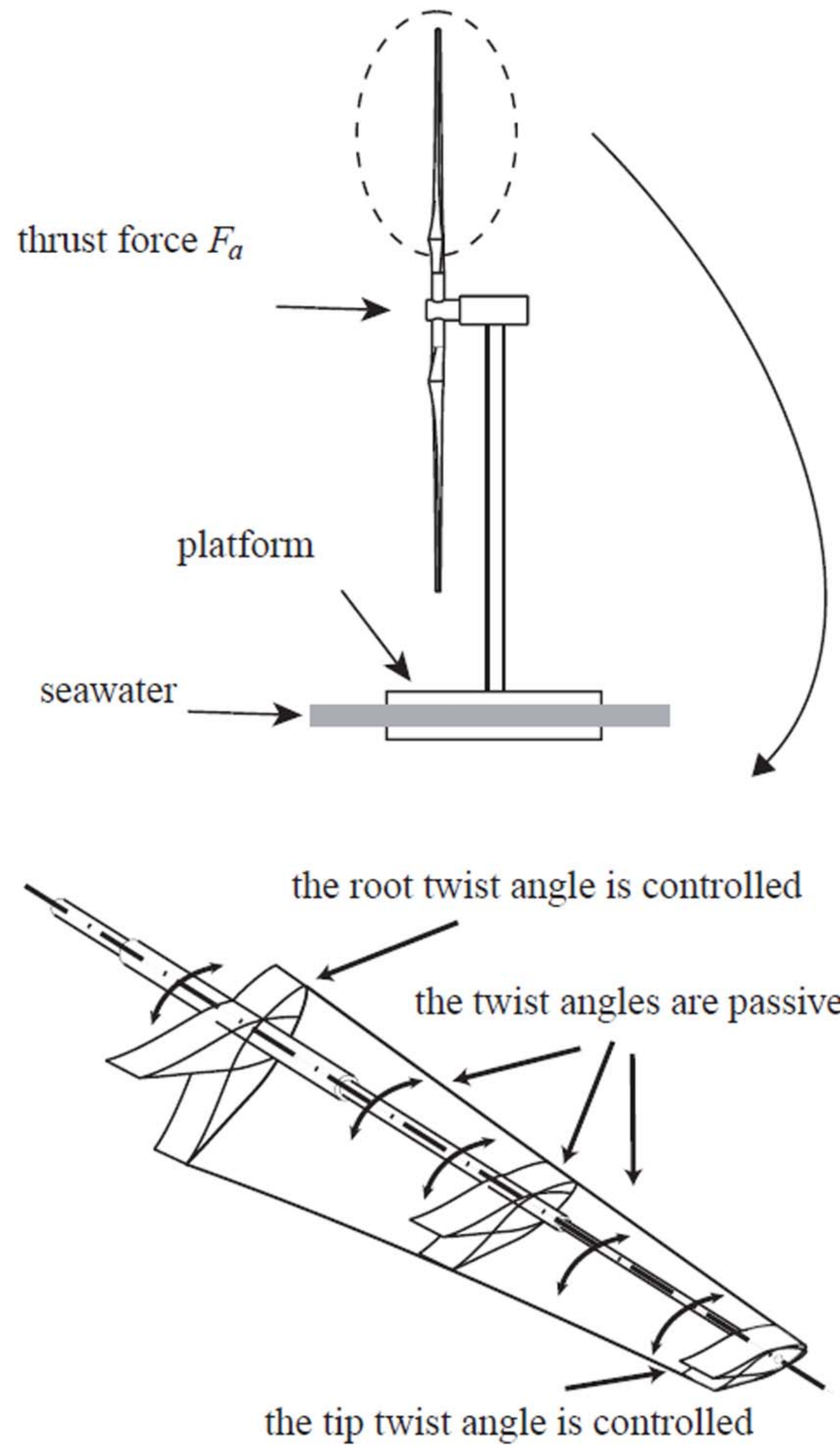}
\caption{Schematic of a floating HAWT rotor with two simplified morphing blades}
\label{wind_turbine}
\end{center}
\end{figure}

The optimization problem at hand has two objective functions: the produced power~$P$ by the wind turbine that should be maximized and the thrust force~$F_a$ that should be minimized. Moreover, the produced power~$P$ should be higher than 1~kW and lower than 25~kW.

$P$ and $F_a$ are calculated based on the Blade Element Momentum Theory (BEMT) knowing the design parameters. A simplified morphing blade with a constant profile type (S809) along the span is used. For more details, the reader is referred to~\cite{Wang2012a}.

To make a good comparison, we take an optimum result as a reference blade, and its parameters are the initial design parameters~\cite{Wang2012a}. The sum of the root chord length~$c_r$ and the tip chord length~$c_t$ is supposed to be constant and equal to~$1.095~\textrm{m}$. Then the four DVs are: ($i$)~the root twist angle~$\gamma_r$; ($ii$)~the tip twist angle~$\gamma_t$; ($iii$)~the chord length at the root~$c_r$; ($iv$)~the rotor rotational speed~$\omega$. The initial values, the lower and upper bounds of the DVs are given in Tab.~\ref{design_Variables}.

\begin{table*}[t]
\caption{The design variables}
\begin{center}
\label{design_Variables}
\small
\begin{tabular}{ c c c c c  }
               \hline
               \hline
  DV & Initial Value & Minimum &  Maximum & Noise values  \\
  \hline
  root twist angle $\gamma_r$ (deg)  &  22.8  & 0 & 35 &   $ \pm 1$ ~( UD) \\
  tip twist angle $\gamma_t$ (deg) &  3.61  & -5 & 15 &   $ \pm 0.5$  ~( UD)\\
  root chord length $c_r$ (m)  &  0.737  & 0.595 & 0.895 &   $ \pm 0.005$ ~( UD) \\
  rotor rotational speed $\omega$ (rpm) & 72 & 40 & 100 &   $ \pm 2$ ~( UD)  \\
               \hline
                \hline
             \end{tabular}
\end{center}
\end{table*}

The other design parameters, such as the number of blades ($b$), tip radius ($r_t$), root radius ($r_r$), air density ($\rho$) and reference wind speed ($v_{re}$) are taken as DEPs. Their nominal values are fixed as shown in Tab.~\ref{design_parameters}.

\begin{table}[t]
\caption{The design environment parameters}
\begin{center}
\label{design_parameters}
\small
\begin{tabular}{ c c c }
               \hline
               \hline
  DEP & Value  & Noise\\
  \hline
  number of blades b & 2  & N/A \\
  tip radius $r_t$ (m) & 5 & $\pm 0.05$ ~( UD) \\
  root radius  $r_r$ (m) & 1.27 & $\pm 0.005$ ~( UD)\\
  air density $\rho$ $(\textrm{kg}/\textrm{m}^3)$ & $1.25$ & $\pm 0.05 $ ~( UD)\\
  reference wind speed  $v_{re}$ (m/s) & 10 & $  \pm 4 $ ~( ND)  \\
               \hline
                \hline
             \end{tabular}
\end{center}
\end{table}

A MOOP is formulated as follows:
\begin{subequations} \label{WT1}
\begin{eqnarray}
  \textrm{minimize } &  f_1({\bf x},{\bf p})  = -P \\
  \textrm{minimize } &  f_2({\bf x},{\bf p})  = F_a\\
  \textrm{over}& {\bf x}=[\gamma_r~\gamma_t~c_r~\omega]  \label{WT-x}\\
~&{\bf p}=[b~r_t~r_r~\rho~v_{re}] \label{WT-p}\\
  \textrm{subject to}  & 1~\textrm{kW} \leq P \leq 25~\textrm{kW} \label{WT-1}\\
~ & 0 ~\textrm{deg} \leq \gamma_r \leq 35~\textrm{deg } \label{WT-hb} \\
~ & -5 ~\textrm{deg} \leq \gamma_t \leq 15~\textrm{deg} \label{WT-lt} \\
~ & 0.595~\textrm{m} \leq c_r \leq 0.895~\textrm{m}   \\
~ & 40~\textrm{rpm }\leq \omega \leq 100~\textrm{rpm} \\
~& b=2;~ r_t=5~\textrm{m};~ r_r = 1.27~\textrm{m}; \\
~& \rho =1.25~\textrm{kg}/\textrm{m}^3; v_{re}=10~\textrm{m/s}
\end{eqnarray}
\end{subequations}

For HAWTs, the variations in DVs and DEPs are unavoidable. Here, we assume that there are noise values in DVs and  DEPs, as shown in Tables~\ref{design_Variables} and~\ref{design_parameters}. Most of them are small variations. However, the noise values in reference wind speed (the performance functions are functionally dependent on it) is quite large, compared with its nominal value. The noise values of small variations are supposed to be Uniform Distribution~(UD). The noise values of the reference wind speed follow a Normal Distribution~(ND) and the corresponding standard deviation is equal to~$2~\textrm{m}/\textrm{s}$.

\subsection{Results Analysis}
\begin{figure}[!htbp]
\begin{center}
\includegraphics[width=0.48\textwidth]{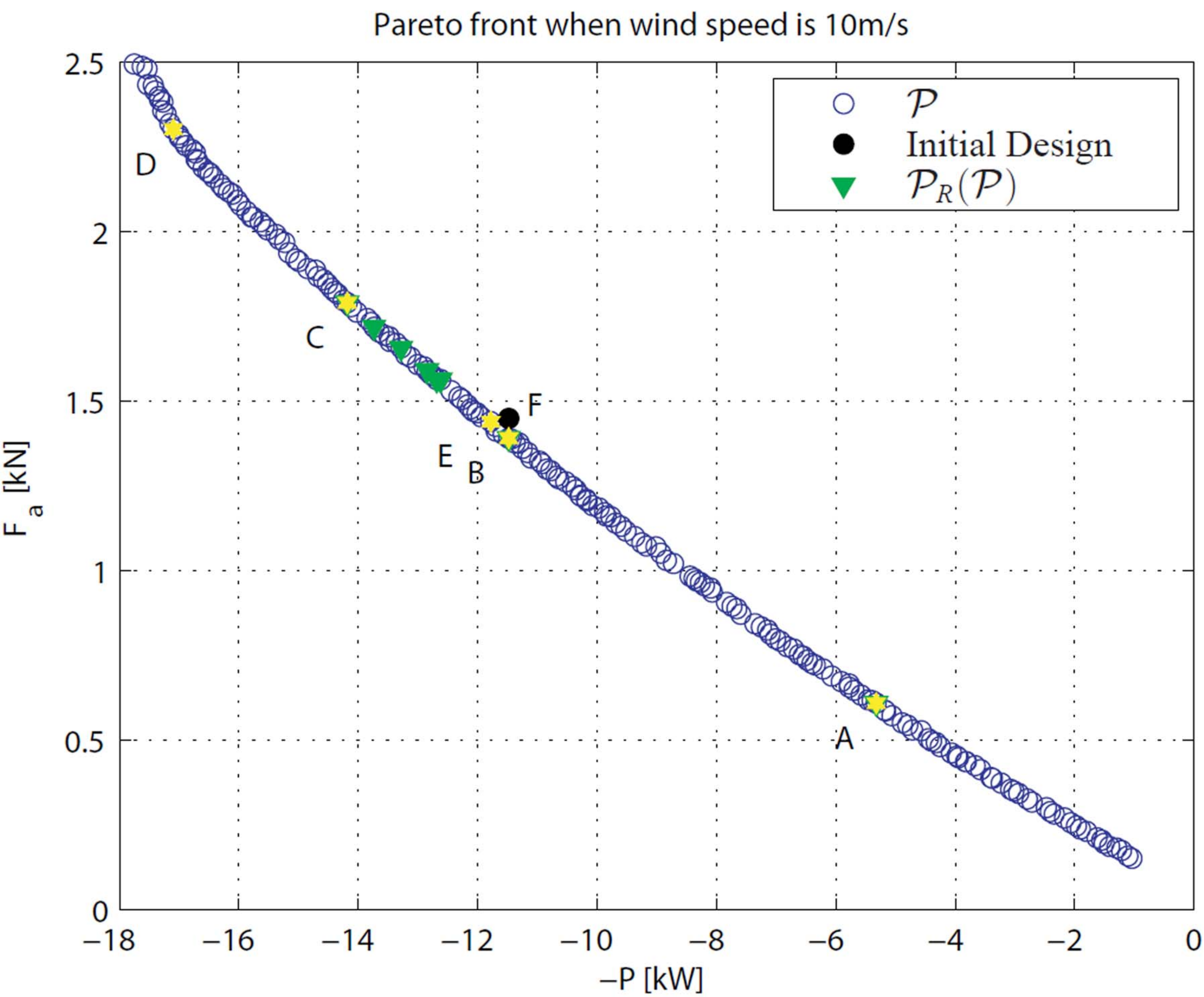}
\caption{The obtained Pareto front for the MOOP of HAWT design}
\label{WT_Morop_PO}
\end{center}
\end{figure}

According to the proposed method, the final solution comes from the Pareto optimal set. In this paper, the Pareto optimal solutions are obtained by using the genetic algorithm (NSGA II ~\cite{Deb2002}). The obtained Pareto front $\mathcal{P}$ for this problem is illustrated in~\Fig {WT_Morop_PO}, including 200~alternative solutions.

To assess the robustness index $I_{RS}$ for each alternative solution, the samples around the nominal values are generated by Latin Hypercube Sampling (LHS) \cite{McKay1979,Luo2008a,Petrone2011}. 1000 samples for each alternative solution are generated with regard to the small variations in DVs and DEPs. Then $I_{RS}$ for each alternative solution can be assessed using with~\Eq {IRS}.

To determine the robustness index $I_{RL}$ with regard to large variations in one DEP for each alternative solution, the new positions of the alternative solutions should be presented, when DEPs change from initial values to the new values. Then the feasibility index and the Pareto optimality Index of a solution can be determined.

In this problem, the wind speed has a continuous probability distribution, varying from 6~m/s to 14~m/s. To simplify the problem, we calculate the probability distribution of the wind speed, and convert it into a discrete probability distribution problem. Table~\ref{table-pdfws} shows the probability distribution of the wind speed. The solution having the maximum PDF amongst the N probabilities is appointed as the initial DEPs ${\bf p}_0$, hence ${\bf p}_0={\bf p}_5$ in this problem.
\begin{table*}[!htbp]
\caption{Probability distributions of wind speed}
\label{table-pdfws}
\begin{center}
\small
\renewcommand{\arraystretch}{2}
\begin{tabular}{ c c c c c c c c c c }
 \hline
  \hline
~& ${\bf p}_1$ & ${\bf p}_2$ &  ${\bf p}_3$ &  ${\bf p}_4$ &  ${\bf p}_5~({\bf p}_0)$ &  ${\bf p}_6$ &  ${\bf p}_7$ &  ${\bf p}_8$ &  ${\bf p}_9$      \\
  \hline
  $v_{re}~[m/s]$ & 6  &   7   &   8  &  9  & 10 &  11  &  12 & 13 &   14   \\
   \hline
  h(\bf p)   & 0.028  &   0.066   & 0.124    &  0.180  & 0.204 &   0.180  &  0.124  & 0.066 &   0.028  \\
  \hline
   \hline
\end{tabular}
\end{center}
\end{table*}

\begin{figure*}[!htbp]
\begin{center}
\includegraphics[width=1\textwidth]{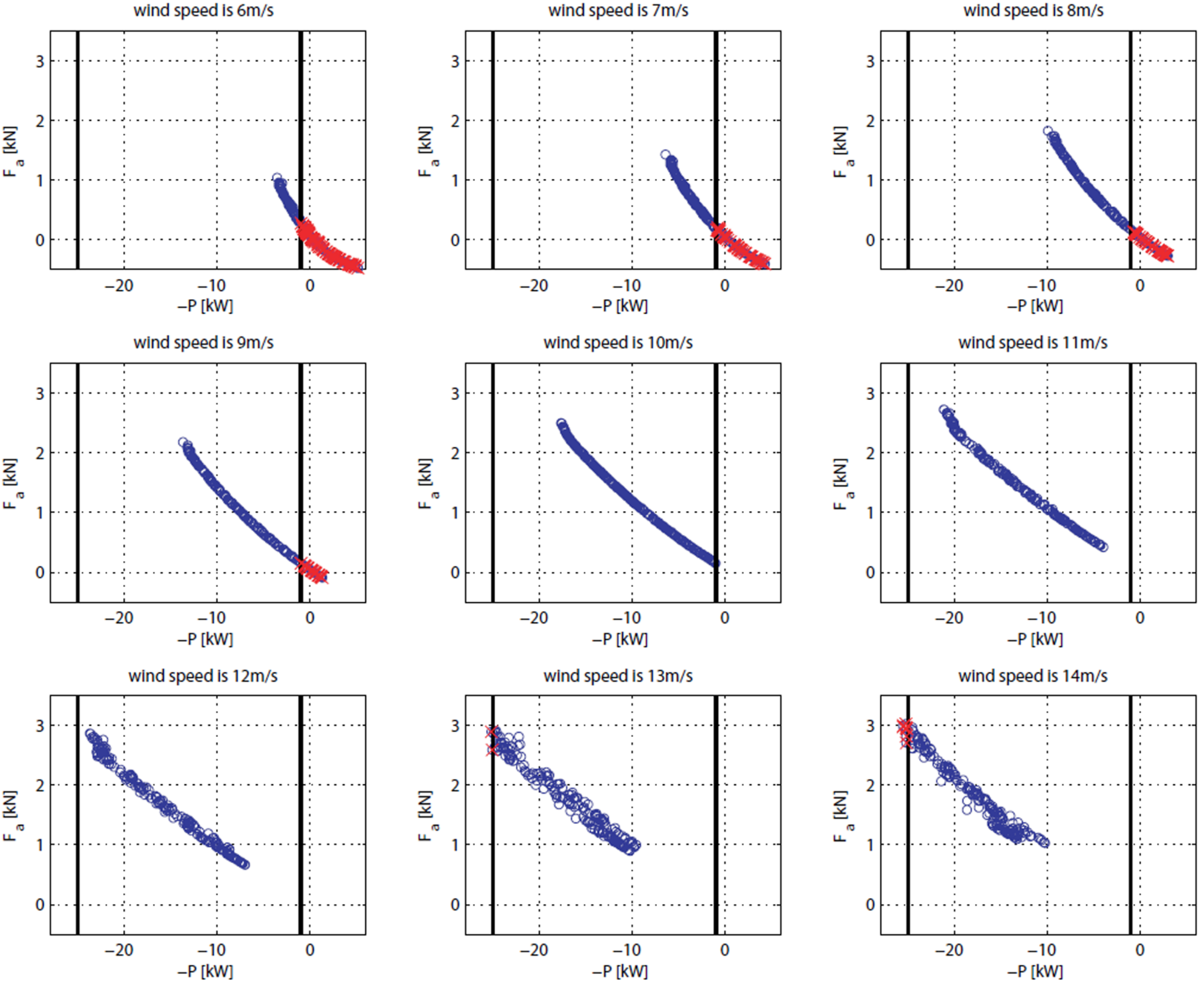}
\caption{The positions of the alternative solutions in different design environments}
\label{WT_Morop_IRL}
\end{center}
\end{figure*}

\Ff {WT_Morop_IRL} illustrates the positions of the alternative solutions in different design environments. The vertical lines represent the constraints in the first objective function, i.e., the produced power $P$, then the solutions between the two lines are feasible. From the results, we can see that some alternative solutions are non-feasible in the new environments. They are illustrated by red crosses in \Fig {WT_Morop_IRL}. The index $I_F$ defined in \Eq {Robustness_indexIF} is null for those solutions. In the new environment, some alternative solutions are dominated by other alternative solutions, then their index $I_P$ defined in \Eq {Robustness_indexIP} is lower than one, but greater than zero. Then the $I_{RL}$ for each alternative solution can be calculated with \Eq {Robustness_indexIRL}.

\begin{figure}[!htbp]
\begin{center}
\scriptsize
\includegraphics[width=0.48\textwidth]{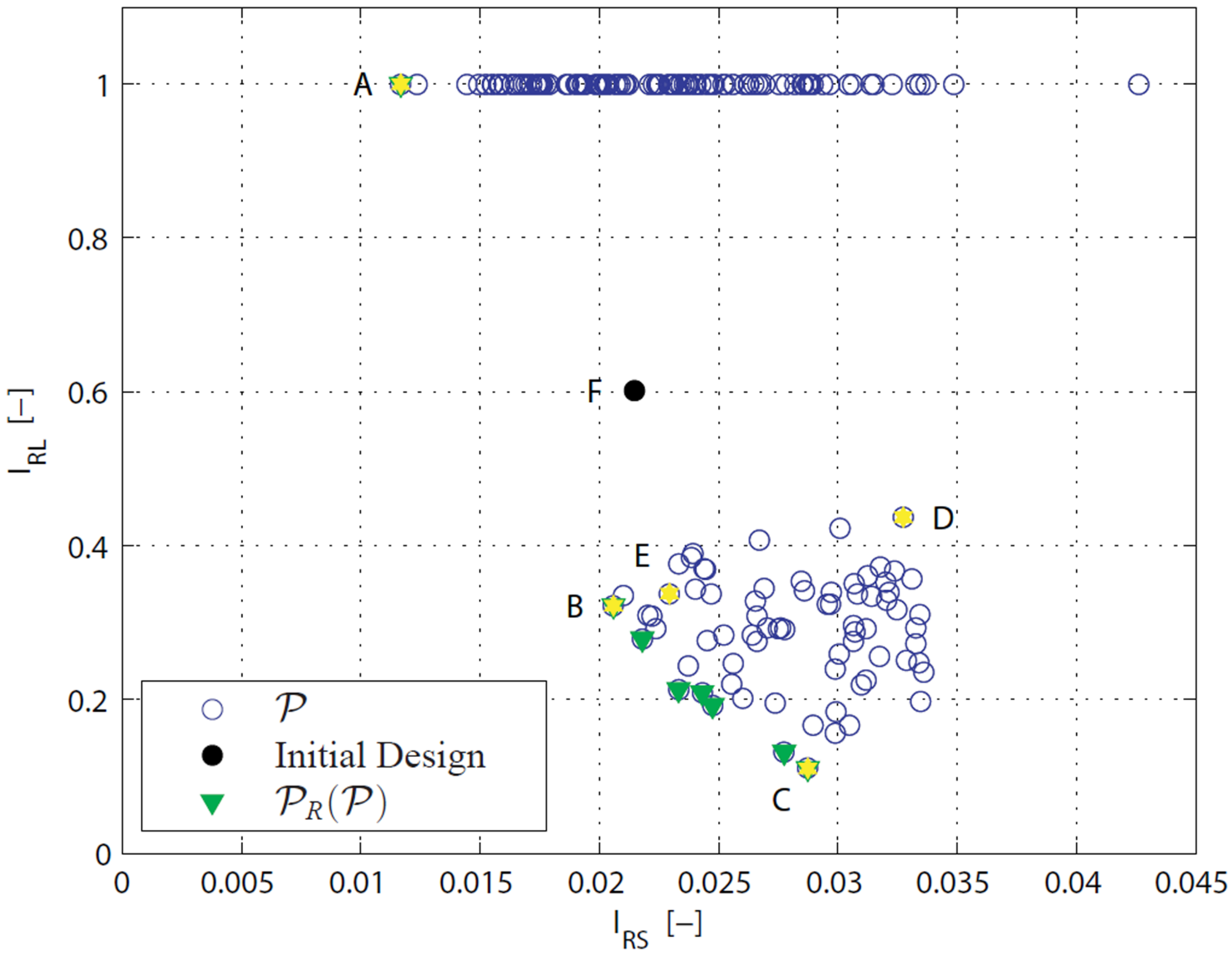}
\caption{The positions of the alternative solutions and initial design in the RF-Space}
\label{WT_Morop_RFSpace}
\end{center}
\end{figure}

\Ff {WT_Morop_RFSpace} illustrates the positions of the alternative solutions in the RF-Space. In this space, one dimension is $I_{RS}$, another dimension is $I_{RL}$. In this example, $I_{RS}$ and $I_{RL}$ are two conflicting objectives, as shown in \Fig {WT_Morop_RFSpace}. A new Pareto front in the RF-Space appears, which represents the most robust solutions amongst the Pareto optimal solutions and denoted by $\mathcal{P}_R(\mathcal{P})$.

\begin{figure}[!htbp]
\begin{center}
\small
\includegraphics[width=0.48\textwidth]{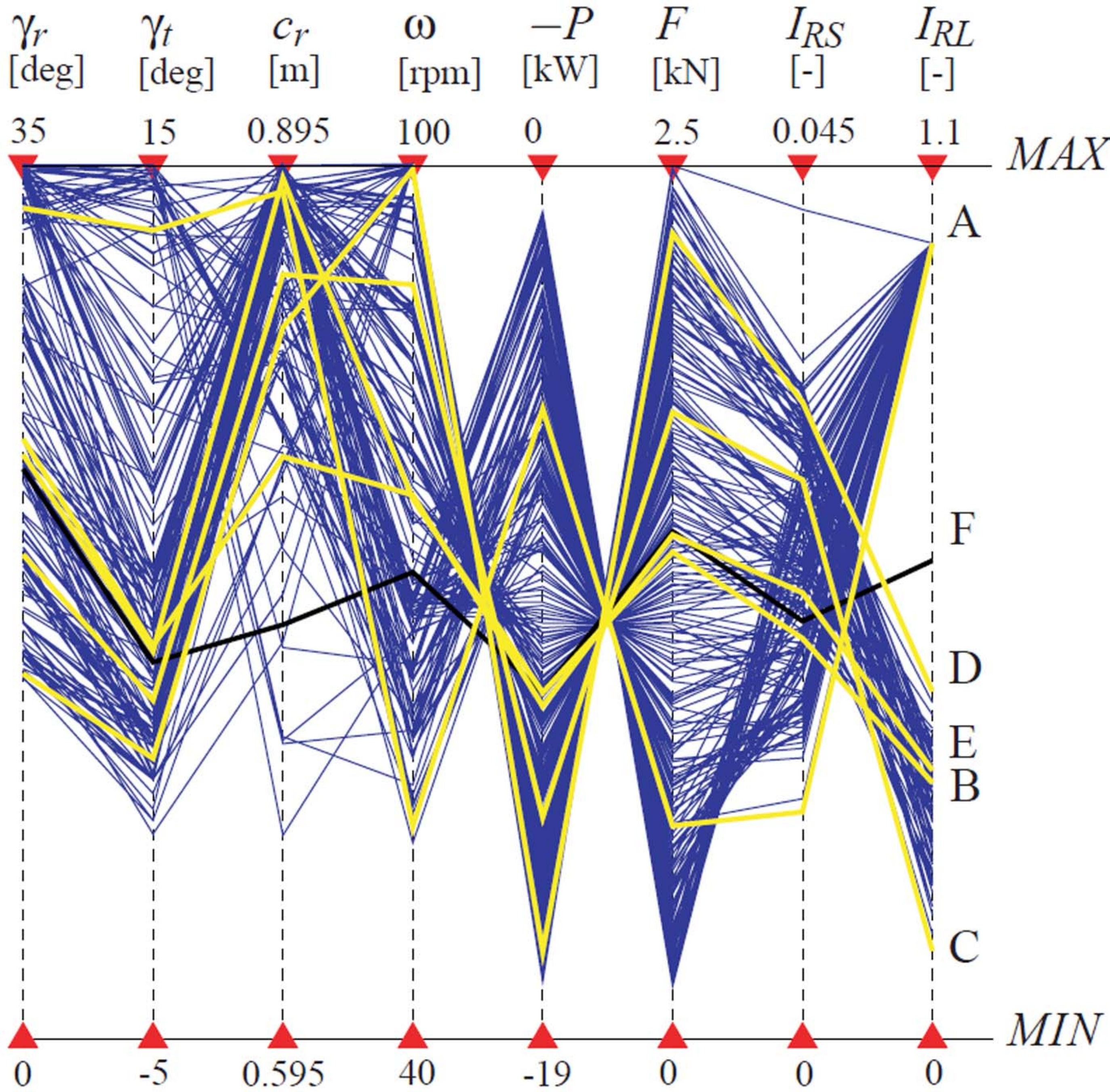}
\caption{Comparison of all alternative solutions and solutions A, B, C, D, E, F in the Decision space, Performance Function space and Robustness Function space}
\label{WT_Morop_Matrix}
\end{center}
\end{figure}

For a good understanding of the two proposed robustness indices defined by \Eq {IRS} and \Eq {Robustness_indexIRL}, the six solutions A, B, C, D, E and~F are selected. Their corresponding positions in the PF-Space and the RF-Space are illustrated in \Fig {WT_Morop_PO} and \Fig {WT_Morop_RFSpace}. All alternative solutions and the initial design are shown in \Fig {WT_Morop_Matrix}. The initial design, solution F, is plotted in black and the other five solutions are plotted in yellow.

\begin{figure*}[!htbp]
\begin{center}
\scriptsize
\includegraphics[width=1\textwidth]{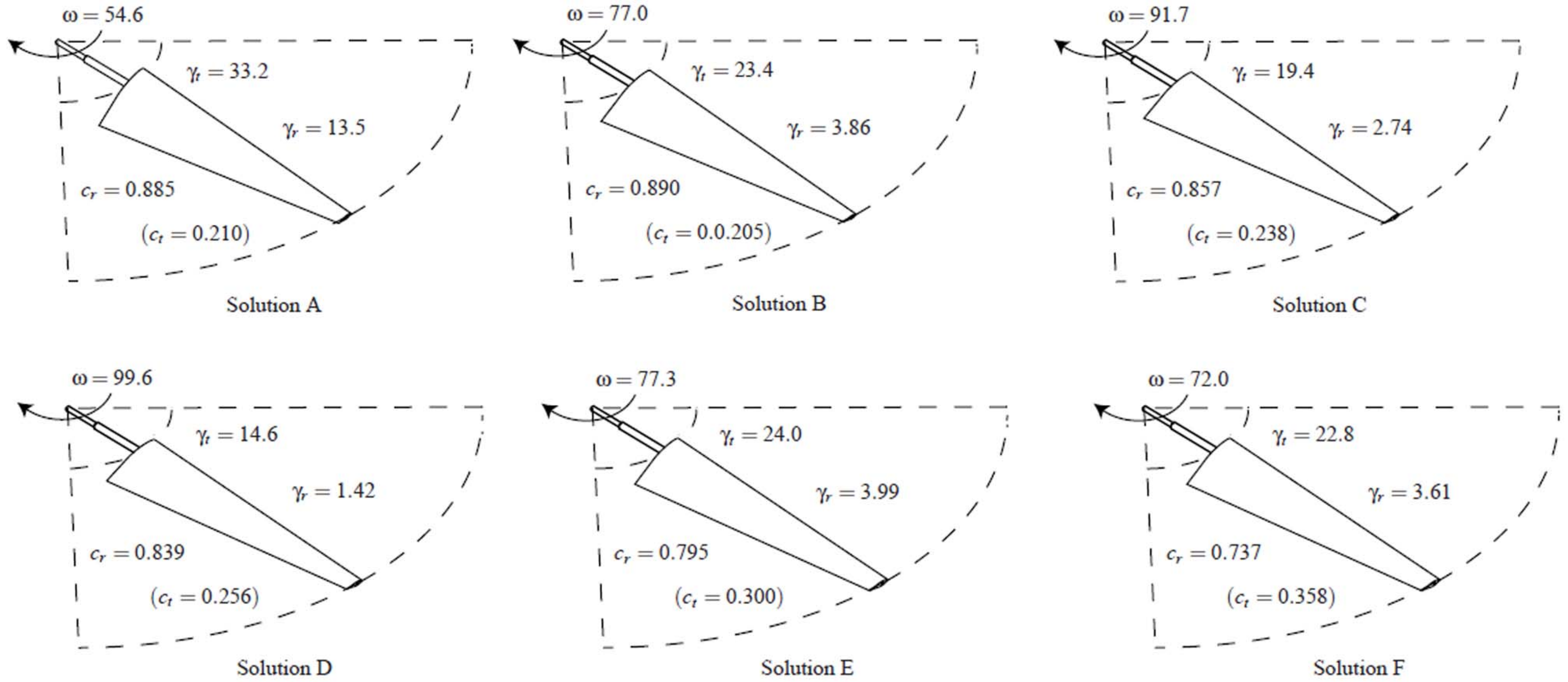}
\caption{3D models for the six selected solutions A, B, C, D, E and~F}
\label{WT_Morop_3Dmodel}
\end{center}
\end{figure*}

Solution A has the minimum $I_{RS}$~value in the $ \mathcal{P}_R(\mathcal{P})$. Solution~B dominates Solution~F in the RF-Space. Solution~C has the minimum $I_{RL}$~value in the $\mathcal{P}$. Solution~D can be dominated by many solutions such as Solutions B, C and E in the RF-Space. Solution~E dominates the initial design in the PF-Space. The 3D models of the six selected solutions and the values of their DVs are depicted in \Fig {WT_Morop_3Dmodel}.

Note that the value of index $I_P$ of a solution depends on its ranking amongst the Pareto optimal solutions. To determine the index $I_P$ of the initial design, i.e., solution~F, we include it into the Pareto optimal set (including 200~alternative solutions) and rank the 201 individuals to generate its ranking in different DEPs. The performance function values and the RIs for these six solutions are given in Tab.~\ref{table6solutions}. For a better understanding of index $I_{RL}$, indices $I_F$ and $I_P$ of these six solutions are also given in Tab.~\ref{table6solutions}.

\begin{table*}[!htbp]
\caption{Comparison of different solutions}
\label{table6solutions}
\begin{center}
\resizebox{0.9\textwidth}{!}{
\renewcommand{\arraystretch}{2}
\begin{tabular}{ c c c c c c c c c c c c c c c}
 \hline
  \hline
Solutions  & -$P$ (kW)  &  $F_a$ (kN) & $I_{RS}$ & $I_F({\bf x})$ & $I_P({\bf x},{\bf p}_1)$ &   $I_P({\bf x},{\bf p}_2)$ &  $I_P({\bf x},{\bf p}_3)$ & $I_P({\bf x},{\bf p}_4)$ &   $I_P({\bf x},{\bf p}_5)$ &  $I_P({\bf x},{\bf p}_6)$ &$I_P({\bf x},{\bf p}_7)$ &   $I_P({\bf x},{\bf p}_8)$ &  $I_P({\bf x},{\bf p}_9)$ & $I_{RL}$  \\
  \hline
  A  & -5.322& 0.609 & 0.0115 & 0  &  -   &   -  &  -  & - &   -  &  -  & - &   -   &  -   & 1  \\
  B  & -11.5 & 1.399 & 0.0208 & 1  &   1   &   1  &  1  & 1 &   1  &  1/4  & 1/6 &  1/10  &  1/4  & 0.322  \\
  C  & -14.2 & 1.799 & 0.0287 & 1  &   1/4   &   1 &  1  & 1/2 &   1  &  1  & 1&   1 &  1  & 0.111 \\
  D  & -17.1 & 2.300 & 0.033  & 1   &   1   &   1  &  1  & 1/2 &   1  &  1/4 & 1/11 &   1/13  &  1/8  & 0.437 \\
  E  & -11.8 & 1.444 & 0.023  & 1   &   1   &   1  &  1/2 & 1/2 & 1 &   1  &  1/7  & 1/9 &   1/3  &  0.338  \\
  F  & -11.5 & 1.455 & 0.0219 & 1   &   1   &   1  &  1/2  & 1 &   1/5  &  1/11  & 1/13 &   1/16  &  1/19  & 0.601 \\
  \hline
   \hline
\end{tabular}}
\end{center}
\end{table*}

With regard to the small variations in DVs and DEPs, each solution has 1000 generated samples in the PF-Space. \Ff {WT_Morop_6solutions_IRs} shows the samples and the nominal values of these six solutions. It is apparent that Solution~A is the most robust one amongst the six solutions considering small variations in DVs and DEPs as the dispersion of the two performance functions is a minimum with regard to small variations in DVs and DEPs for this solution. It matches with their positions in the RF-Space. The robustness indices $I_{RS}$ with regard to small variations for the five other solutions are given in the fourth column of Tab.~\ref{table6solutions}. Here is the order of the six solutions from the most robust one to the least robust one with respect to $I_{RS}$: (1)~A; (2)~B; (3)~F; (4)~E; (5)~C; (6)~D.

\begin{table}[!htbp]
\caption{Comparison of six solutions with regard to variations in DV and DEPs}
\label{table5solutions_WTeg}
\begin{center}
\small
\resizebox{0.5\textwidth}{!}{
\renewcommand{\arraystretch}{2}
\begin{tabular}{ c c c c c c c}
 \hline
  \hline
Solutions  & $-P$ (kW)  &  $F_a$ (kN) &  $\Delta P^S$  &  $\Delta F_a^S $  & $\Delta  P^L$   &  $\Delta F_a^L$  \\
  \hline
A&   -5.32&    0.609&    0.576& 0.0694&    6.28&   0.416\\
B&   -11.5&    1.399&    0.801&   0.110&    10.2&   0.410\\
C&   -14.2&    1.799&    1.22&    0.168&    12.5&    0.557\\
D&   -17.1&    2.300&    1.08&    0.177&    13.9&    0.493\\
E&   -11.8&    1.444&    0.977&   0.120&    10.7&   0.453\\
F&   -11.5&    1.455&    0.910&    0.113&    9.57&   0.366\\
  \hline
  \hline
\end{tabular}}
\end{center}
\end{table}

With regard to large variations in DEPs, all Pareto optimal solutions have sufficiently large variations in performance functions. In Tab.~\ref{table5solutions_WTeg}, $-P$ and $F_a$ represent the nominal performance function values of the alternative solutions.  $\Delta P^S$ and $\Delta F_a^S$ represent the largest distance of the actual performance function values (1000 samples for each alternative solution) and nominal performance function values in $P$ and $F_a$ respectively, with regard to the small variations in DVs and DEPs. $\Delta P^L$ and $\Delta F_a^L$ represent the largest distance of the actual performance function values (9 samples for each alternative solution) and nominal performance function values in $P$ and $F_a$ respectively, with regard to large variations in DEPs.
The results show that $\Delta P^L>> \Delta P^S$ and $\Delta F_a^L>> \Delta F_a^S$, so with regard to large variations in DEPs, the traditional methods (results based on the difference of actual performance function values and nominal ones) make little sense. On the contrary, the method proposed in this paper is relevant. Indeed, here is the order of the six solutions from the most robust one to the least robust one based on the $I_{RL}$ index: (1)~C; (1)~B; (3)~E; (4)~D; (5)~F; (6)~A.

\begin{figure*}[!htbp]
\begin{center}
\includegraphics[width=1\textwidth]{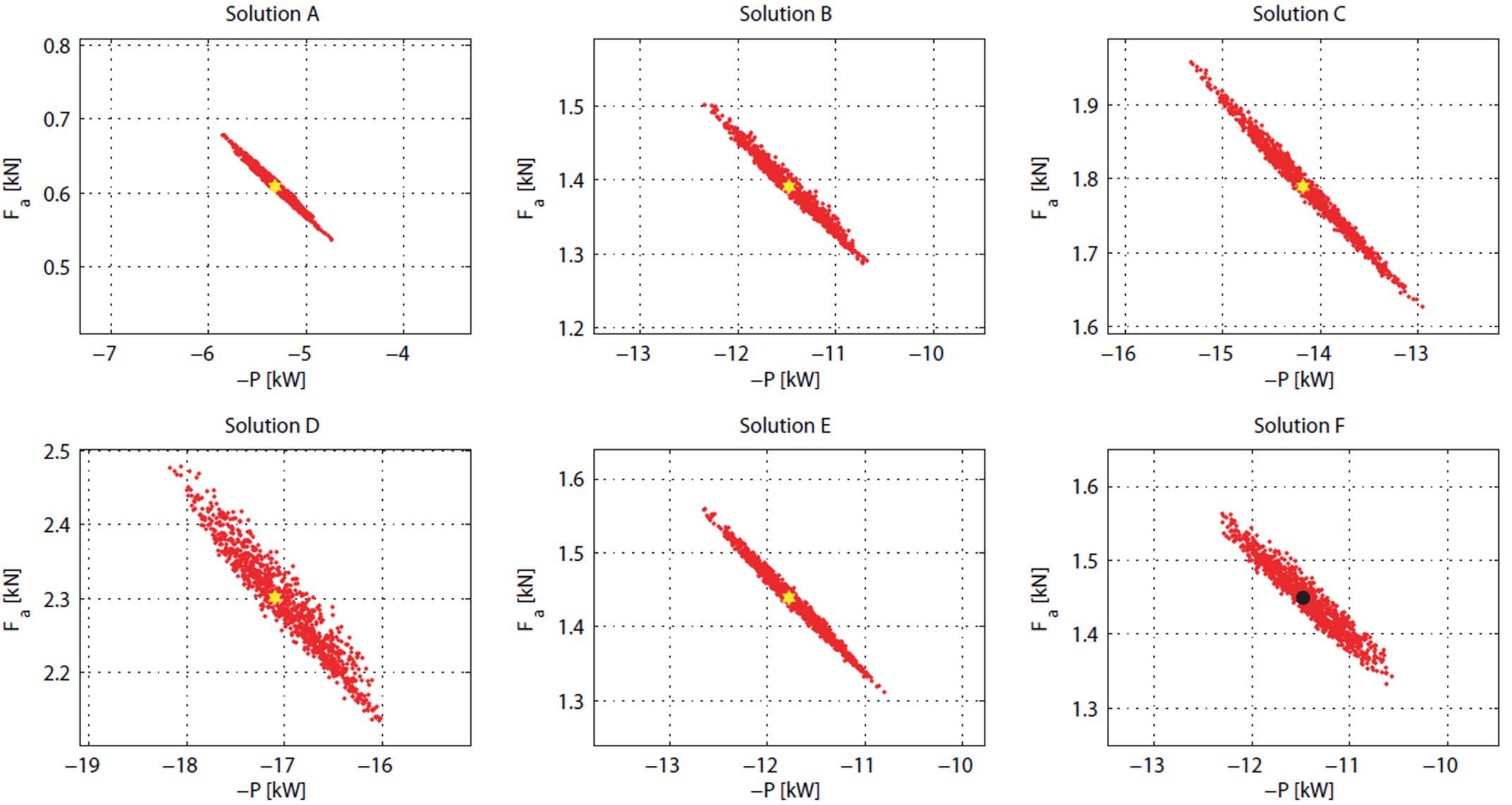}
\caption{The generated samples and the nominal values of these six solutions with regard to small variations in DVs and DEPs}
\label{WT_Morop_6solutions_IRs}
\end{center}
\end{figure*}

In summary, the designer can select the final solution from the new Pareto front in the RF-Space, according to his/her requirement. For example, if the designer prefers the $I_{RS}$, solution~A should be selected. If the designer prefers the $I_{RL}$, solution~C should be selected.

\section{Conclusions and Future Work}

In this paper, a new method for solving Multi-Objective Robust Optimization Problems~(MOROP) has been introduced and two illustrative examples have been given to highlight the main contributions of the paper. Two Robustness Indices~(RI) have been introduced to deal with MOROP where not only small variations in Design Variables (DVs) and Design Environment Parameters~(DEPs) are considered, but also large variations in DEPs. The first robustness index, named~$I_{RS}$, characterizes the robustness of MOOP against small variations in DVs and DEPs. The second robustness index, named $I_{RL}$, characterizes the robustness of Multi-Objective Optimization Problems~(MOOP) against large variations in DEPs. The robustness index~$I_{RS}$ is calculated based on the standard deviations and the differences between the expected values and nominal values of the performances. The smaller $I_{RS}$, the more robust the design. The robustness index~$I_{RL}$ is calculated based on the solution's ability to be optimal in different design environments. The smaller $I_{RL}$, the more robust the design. To make a trade-off between the two proposed robustness indices, a concept of Robust Function Space (RF-Space) has been introduced. Then each Pareto optimal solution has a position in the RF-Space. The designer can select the final solution from the Pareto optimal set based on its new position in the RF-Space. However, some problems are still not solved. For instance, it is not always an easy task to claim whether a variation is small or large. Moreover, new formulations for $I_{RS}$ and $I_{RL}$ should be discussed in the future.



\begin{thebibliography}{10}

\bibitem{Caro2005}
Caro, S., Bennis, F., and Wenger, P., 2005.
\newblock ``{Tolerance synthesis of mechanisms: a robust design approach}''.
\newblock {\em Journal of Mechanical Design, {\bf 127}}, pp.~86--94.

\bibitem{Beyer2007}
Beyer, H., and Sendhoff, B., 2007.
\newblock ``{Robust optimization-a comprehensive survey}''.
\newblock {\em Computer methods in applied mechanics and engineering, {\bf
  196}}(33-34), July, pp.~3190--3218.

\bibitem{Deb2005}
Deb, K., and Gupta, H., 2005.
\newblock ``{Searching for robust Pareto-optimal solutions in multi-objective
  optimization}''.
\newblock In {\em Evolutionary Multi-Criterion Optimization}, C.~{Coello
  Coello}, A.~{Hern\'{a}ndez Aguirre}, and E.~Zitzler, eds. Springer Berlin /
  Heidelberg, pp.~150--164.

\bibitem{Deb2006}
Deb, K., and Gupta, H., 2006.
\newblock ``{Introducing robustness in multi-objective optimization}''.
\newblock {\em Evolutionary Computation, {\bf 14}}(4), pp.~463--494.

\bibitem{Saha2011}
Saha, A., and Ray, T., 2011.
\newblock ``{Practical Robust Design Optimization Using Evolutionary
  Algorithms}''.
\newblock {\em Journal of Mechanical Design, {\bf 133}}(10), p.~101012.

\bibitem{Gaspar-Cunha2006}
Gaspar-Cunha, A., and Covas, J., 2006.
\newblock ``{Robustness using Multi-Objective Evolutionary Algorithms}''.
\newblock In {\em Applications of Soft Computing}, K.~{Tiwari, Ashutosh and
  Roy, Rajkumar and Knowles, Joshua and Avineri, Erel and Dahal}, ed. Springer
  Berlin / Heidelberg, pp.~353--362.

\bibitem{Gaspar-Cunha2008}
Gaspar-Cunha, A., and Covas, J., 2008.
\newblock ``{Robustness in multi-objective optimization using evolutionary
  algorithms}''.
\newblock {\em Computational Optimization and Applications, {\bf 39}}(1), June,
  pp.~75--96.

\bibitem{Barrico2006}
Barrico, C., and Antunes, C., 2006.
\newblock ``{A New Approach to Robustness Analysis in Multi-Objective
  Optimization}''.
\newblock In 7th International Conference on Multi-Objective Programming and
  Goal Programming (MOPGP 2006), Loire Valley, City of Tours, France, no.~x,
  pp.~12--15.

\bibitem{Barrico2006a}
Barrico, C., and Antunes, C.~H., 2006.
\newblock ``{Robustness analysis in multi-obejective optimization Using a
  degree of robustness concept}''.
\newblock In IEEE Congress on Evolutionary Computation, pp.~1887--1892.

\bibitem{Giassi2004}
Giassi, A., Bennis, F., and Maisonneuve, J.-J., 2004.
\newblock ``{Multidisciplinary design optimisation and robust design approaches
  applied to concurrent design}''.
\newblock {\em Structural and Multidisciplinary Optimization, {\bf 28}}(5),
  Aug., pp.~356--371.

\bibitem{Caro2005a}
Caro, S., Bennis, F., and Wenger, P., 2005.
\newblock ``{Comparison of Robustness Indices and Introduction of a Tolerance
  Synthesis Method for Mechanisms}''.
\newblock In Canadian Congress of Applied Mechanics, pp.~3--4.

\bibitem{Wang2012}
Wang, W., Caro, S., Bennis, F., and Augusto, O.~B., 2012.
\newblock ``{Toward The Use Of Pareto Performance Solutions And Pareto
  Robustness Solutions For Multi-Objective Robust Optimization Problems}''.
\newblock In Proceedings of the ASME 2012 11th Biennial Conference On
  Engineering Systems Design And Analysis, pp.~ESDA2012--82099.

\bibitem{Augusto2012a}
Augusto, O., Bennis, F., and Caro, S., 2012.
\newblock ``{Multiobjective engineering design optimization problems: a
  sensitivity analysis approach}''.
\newblock {\em Pesquisa Operacional}, pp.~575--596.

\bibitem{Gunawan2005}
Gunawan, S., and Azarm, S., 2005.
\newblock ``{Multi-objective robust optimization using a sensitivity region
  concept}''.
\newblock {\em Structural and Multidisciplinary Optimization, {\bf 29}}(1),
  Aug., pp.~50--60.

\bibitem{Li2005}
Li, M., Azarm, S., and Aute, V., 2005.
\newblock ``{A multi-objective genetic algorithm for robust design
  optimization}''.
\newblock In Proceedings of the 2005 conference on Genetic and evolutionary
  computation, June, pp.~25--29.

\bibitem{Li2006}
Li, M., Azarm, S., and Boyars, A., 2006.
\newblock ``{A new deterministic approach using sensitivity region measures for
  multi-objective robust and feasibility robust design optimization}''.
\newblock {\em Journal of mechanical design, {\bf 128}}(4), pp.~874--883.

\bibitem{Hu2009}
Hu, W., Li, M., Azarm, S., {Al Hashimi}, S., Almansoori, A., and Al-Qasas, N.,
  2009.
\newblock ``{Improving Multi-Objective Robust Optimization Under Interval
  Uncertainty Using Worst Possible Point Constraint Cuts}''.
\newblock {\em ASME Conference Proceedings, {\bf 2009}}, pp.~1193--1203.

\bibitem{Hung2011}
Hung, T., and Chan, K., 2011.
\newblock ``{Multi-objective design and tolerance allocation for single-and
  multi-level systems}''.
\newblock {\em Journal of Intelligent Manufacturing}(0956-5515), Nov.,
  pp.~1--15.

\bibitem{Hung2013}
Hung, T., and Chan, K., 2013.
\newblock ``{Uncertainty quantifications of Pareto optima in multiobjective
  problems}''.
\newblock {\em Journal of Intelligent Manufacturing, {\bf 24}}(0956-5515),
  Nov., pp.~385--395.

\bibitem{Apley2006}
Apley, D.~W., Liu, J., and Chen, W., 2006.
\newblock ``{Understanding the Effects of Model Uncertainty in Robust Design
  With Computer Experiments}''.
\newblock {\em Journal of Mechanical Design, {\bf 128}}(4), pp.~945--958.

\bibitem{Hassan2008}
Hassan, G., and Clack, C., 2008.
\newblock ``{Multiobjective robustness for portfolio optimization in volatile
  environments}''.
\newblock In GECCO '08: Proceedings of the 10th annual conference on Genetic
  and evolutionary computation, ACM Press, pp.~1507--1514.

\bibitem{Hassan2009}
Hassan, G., and Clack, C.~D., 2009.
\newblock ``{Robustness of Multiple Objective GP Stock-Picking in Unstable
  Financial Markets}''.
\newblock In GECCO '09: Proceedings of the 11th Annual conference on Genetic
  and evolutionary computation, pp.~1513----1520.

\bibitem{Lu2009}
Lu, X., and Li, H.-X., 2009.
\newblock ``{Perturbation Theory Based Robust Design Under Model
  Uncertainty}''.
\newblock {\em Journal of Mechanical Design, {\bf 131}}(11), p.~111006.

\bibitem{Moshaiov2006}
Moshaiov, A., and Avigad, G., 2006.
\newblock ``{Concept-based IEC for Multi-objective Search with Robustness to
  Human Preference Uncertainty}''.
\newblock In 2006 IEEE International Conference on Evolutionary Computation,
  Ieee, pp.~1893--1900.

\bibitem{Avigad2011}
Avigad, G., Eisenstadt, E., and Schuetze, O., 2011.
\newblock ``{Handling changes of performance requirements in multi-objective
  problems}''.
\newblock {\em Journal of Engineering Design}(May 2012), Nov., pp.~1--21.

\bibitem{Bui2012}
Bui, L.~T., Abbass, H.~A., and Barlow, M., 2012.
\newblock ``{Robustness Against the Decision-Maker's Attitude to Risk in
  Problems With Conflicting Objectives}''.
\newblock {\em IEEE Transactions on Evolutionary Computation, {\bf 16}}(99),
  pp.~1--19.

\bibitem{Cruz2011}
Cruz, C., Gonz\'{a}lez, J.~R., and Pelta, D.~a., 2011.
\newblock ``{Optimization in dynamic environments: a survey on problems,
  methods and measures}''.
\newblock {\em Soft Computing, {\bf 15}}(7), Dec., pp.~1427--1448.

\bibitem{Farina2003}
Farina, M., Deb, K., and Amato, P., 2003.
\newblock ``{Dynamic multiobjective optimization problems: test cases,
  approximations, and applications}''.
\newblock In {\em Evolutionary Multi-Criterion Optimization}, K.~{Fonseca,
  CarlosM. and Fleming, PeterJ. and Zitzler, Eckart and Thiele, Lothar and
  Deb}, ed., no.~Mi in Lecture Notes in Computer Science. Springer Berlin
  Heidelberg, pp.~311--326.

\bibitem{Deb2007}
Deb, K., N., U. B.~R., and Karthik, S., 2007.
\newblock ``{Dynamic multi-objective optimization and decision-making using
  modified NSGA-II: a case study on hydro-thermal power scheduling}''.
\newblock In {\em Evolutionary Multi-Criterion Optimization}, T.~{Obayashi,
  Shigeru and Deb, Kalyanmoy and Poloni, Carlo and Hiroyasu, Tomoyuki and
  Murata}, ed. Springer Berlin Heidelberg, pp.~803--817.

\bibitem{Chen2011}
Chen, R., and Zeng, W., 2011.
\newblock ``{Multi-objective optimization in dynamic environment: A review}''.
\newblock In 2011 6th International Conference on Computer Science \& Education
  (ICCSE), Vol.~0, Ieee, pp.~78--82.

\bibitem{Nazari2013}
Nazari, V., and Notash, L., 2013.
\newblock ``{Workspace of wire-actuated parallel manipulators and variations in
  design parameters}''.
\newblock {\em Transactions of the Canadian Society for Mechanical Engineering,
  {\bf 37}}(2), pp.~215--229.

\bibitem{Saltelli2008}
Saltelli, A., Ratto, M., Andres, T., and Campolongo, F., 2008.
\newblock {\em {Global sensitivity analysis: the primer}}.
\newblock Wiley-Interscience.

\bibitem{Tomlin2011}
Tomlin, A., and Ziehn, T., 2011.
\newblock ``{The Use of Global Sensitivity Methods for the Analysis, Evaluation
  and Improvement of Complex Modelling Systems}''.
\newblock In {\em Coping with Complexity: Model Reduction and Data Analysis},
  Vol.~75. Springer Berlin Heidelberg, pp.~9--36.

\bibitem{Tong2008}
Tong, C., and Graziani, F., 2008.
\newblock ``{A Practical Global Sensitivity Analysis Methodology for
  Multi-Physics Applications}''.
\newblock In {\em Computational Methods in Transport: Verification and
  Validation}, Vol.~62. pp.~277--299.

\bibitem{Augustin2011}
Augustin, J., 2011.
\newblock ``{Global sensitivity analysis applied to food safety risk
  assessment}''.
\newblock {\em Electronic Journal of Applied Statistical Analysis, {\bf 4}}(2),
  pp.~255--264.

\bibitem{Gunawan2004}
Gunawan, S., and Azarm, S., 2004.
\newblock ``{Non-gradient based parameter sensitivity estimation for single
  objective robust design optimization}''.
\newblock {\em Journal of mechanical design, {\bf 126}}(3), pp.~395--402.

\bibitem{Lu2010}
Lu, X., Li, H.-X., and Chen, C. L.~P., 2010.
\newblock ``{Variable Sensitivity-Based Deterministic Robust Design for
  Nonlinear System}''.
\newblock {\em Journal of Mechanical Design, {\bf 132}}(6), p.~064502.

\bibitem{Deb2001book}
Deb, K., 2001.
\newblock {\em {Multi-objective Optimization using Evolutionary Algorithms}}.
\newblock John Wiley and Sons, New York.

\bibitem{Li.X2009}
Li, X., and Wong, H.-S., 2009.
\newblock ``{Logic optimality for multi-objective optimization}''.
\newblock {\em Applied Mathematics and Computation, {\bf 215}}(8), Dec.,
  pp.~3045--3056.

\bibitem{Coello2006}
{Coello Coello}, C., 2006.
\newblock ``{Evolutionary multi-objective optimization: a historical view of
  the field}''.
\newblock {\em Computational Intelligence Magazine, IEEE}(February 2006),
  pp.~28--36.

\bibitem{Baker2008}
Baker, J.~W., Schubert, M., and Faber, M.~H., 2008.
\newblock ``{On the assessment of robustness}''.
\newblock {\em Structural Safety, {\bf 30}}(3), May, pp.~253--267.

\bibitem{Fink2009}
Fink, G., Steiger, R., and K\"{o}hler, J., 2009.
\newblock ``{Definition of Robustness and related terms}''.
\newblock In Joint Workshop of COST Actions TU0601 and E55, pp.~21--25.

\bibitem{Wieland2012}
Wieland, A., and Wallenburg, C., 2012.
\newblock ``{Dealing with supply chain risks: Linking risk management practices
  and strategies to performance}''.
\newblock {\em International Journal of Physical Distribution \& Logistics
  Management, {\bf 42}}(10).

\bibitem{Kitano2004}
Kitano, H., 2004.
\newblock ``{Biological robustness.}''.
\newblock {\em Nature reviews. Genetics, {\bf 5}}(11), Nov., pp.~826--37.

\bibitem{Felix2008}
F\'{e}lix, M.-a., and Wagner, A., 2008.
\newblock ``{Robustness and evolution: concepts, insights and challenges from a
  developmental model system.}''.
\newblock {\em Heredity, {\bf 100}}(2), Feb., pp.~132--40.

\bibitem{Kuorikoski2007}
Kuorikoski, J., Lehtinen, A., and Marchionni, C., 2007.
\newblock ``{Economics as robustness analysis}''.
\newblock In Models and Simulations 2, pp.~1--29.

\bibitem{Shahrokni2013}
Shahrokni, A., and Feldt, R., 2013.
\newblock ``{A systematic review of software robustness}''.
\newblock {\em Information and Software Technology, {\bf 55}}(1), Jan.,
  pp.~1--17.

\bibitem{Taguchi1993}
Taguchi, G., 1993.
\newblock {\em {Taguchi on robust technology development: bringing quality
  engineering upstream}}.
\newblock ASME Press.

\bibitem{PARK2006}
Park, G.-J., Lee, K.-H., Kwon, H., and Hwang, K., 2006.
\newblock ``{Robust design: an overview}''.
\newblock {\em AIAA journal, {\bf 44}}(1), pp.~181--191.

\bibitem{Arvidsson2008}
Arvidsson, M., 2008.
\newblock ``{Principles of robust design methodology}''.
\newblock {\em Quality and Reliability Engineering International, {\bf 24}}(1),
  pp.~23--35.

\bibitem{Suh2001}
Suh, N., 2001.
\newblock {\em {Axiomatic design: advances and applications}}.
\newblock Oxford Univ. Press, New York.

\bibitem{Box1993}
Box, G.~E., and Fung, C., 1993.
\newblock ``{Is Your Robust Design Procedure Robust?}''.
\newblock {\em Quality Engineering, {\bf 6}}, pp.~503--514.

\bibitem{Zang2005a}
Zang, C., Friswell, M., and Mottershead, J., 2005.
\newblock ``{A review of robust optimal design and its application in
  dynamics}''.
\newblock {\em Computers \& Structures, {\bf 83}}(4-5), Jan., pp.~315--326.

\bibitem{Huang2006}
Huang, B., and Du, X., 2007.
\newblock ``{Analytical robustness assessment for robust design}''.
\newblock {\em Structural and Multidisciplinary Optimization, {\bf 34}}(2),
  Dec., pp.~123--137.

\bibitem{Huang2006a}
Huang, B., and Du, X., 2006.
\newblock ``{A robust design method using variable transformation and
  Gauss–Hermite integration}''.
\newblock {\em International Journal for Numerical Methods in Engineering, {\bf
  66}}(12), June, pp.~1841--1858.

\bibitem{Lee2010}
Lee, M.~C., Mikulik, Z., Kelly, D.~W., Thomson, R.~S., and Degenhardt, R.,
  2010.
\newblock ``{Robust design - A concept for imperfection insensitive composite
  structures}''.
\newblock {\em Composite Structures, {\bf 92}}(6), May, pp.~1469--1477.

\bibitem{Ferreira2008}
Ferreira, J., Fonseca, C.~M., Covas, J.~A., and Gaspar-Cunha, A., 2008.
\newblock ``{Evolutionary Multi-Objective Robust Optimization}''.
\newblock In {\em Advances in Evolutionary Algorithms}, Z.~Xiong, ed.,
  no.~November. InTech, ch.~13.

\bibitem{Barrico2007}
Barrico, C., and Antunes, C., 2007.
\newblock ``{An Evolutionary Approach for Assessing the Degree of Robustness of
  Solutions to Multi-Objective Models}''.
\newblock In {\em Evolutionary Computation in Dynamic and Uncertain
  Environments－Studies in Computational Intelligence}, S.~Yang, Y.-S. Ong,
  and Y.~Jin, eds., Vol.~51 of {\em Studies in Computational Intelligence}.
  Springer Berlin / Heidelberg, pp.~565--582.

\bibitem{Barrico2009}
Barrico, C., Antunes, C., and Pires, D., 2009.
\newblock ``{Robustness Analysis in Evolutionary Multi-Objective Optimization
  Applied to VAR Planning in Electrical Distribution Networks}''.
\newblock In {\em Evolutionary Computation in Combinatorial Optimization},
  C.~Cotta and P.~Cowling, eds., Vol.~5482 of {\em Lecture Notes in Computer
  Science}. Springer Berlin / Heidelberg, pp.~216--227.

\bibitem{Hu2011b}
Hu, W., Li, M., Azarm, S., and Almansoori, A., 2011.
\newblock ``{Multi-Objective Robust Optimization Under Interval Uncertainty
  Using Online Approximation and Constraint Cuts}''.
\newblock {\em Journal of Mechanical Design, {\bf 133}}(6), p.~061002.

\bibitem{Nejlaoui2011}
Nejlaoui, M., Houidi, A., Affi, Z., and Romdhane, L., 2011.
\newblock ``{Multiobjective robust design optimization of rail vehicle}''.
\newblock In 13th World Congress in Mechanism and Machine Science, p.~A15\_412.

\bibitem{Nejlaoui2012}
Nejlaoui, M., Houidi, A., Affi, Z., and Romdhane, L., 2012.
\newblock ``{Multiobjective robust design optimization of rail vehicle moving
  in short radius curved tracks based on the safety and comfort criteria}''.
\newblock {\em Simulation Modelling Practice and Theory}, Oct.

\bibitem{Erfani2010}
Erfani, T., and Utyuzhnikov, S., 2010.
\newblock ``{Handling Uncertainty and Finding Robust Pareto Frontier in
  Multiobjective Optimization Using Fuzzy Set Theory}''.
\newblock In 51st AIAA/ASME/ASCE/AHS/ASC Structures, Structural Dynamics, and
  Materials Conference et al, no.~April, pp.~1--7.

\bibitem{Erfani2012}
Erfani, T., and Utyuzhnikov, S., 2012.
\newblock ``{Control of robust design in multiobjective optimization under
  uncertainties}''.
\newblock {\em Structural and Multidisciplinary Optimization, {\bf
  45}}(1615-147X), Aug., pp.~247----256.

\bibitem{AitBrik2006}
{Ait Brik}, B., Ghanmi, S., Bouhaddi, N., and Cogan, S., 2006.
\newblock ``{Robust Design in Structural Mechanics}''.
\newblock {\em International Journal for Computational Methods in Engineering
  Science and Mechanics, {\bf 8}}(1), Dec., pp.~39--49.

\bibitem{Ghanmi2007}
Ghanmi, S., Bouazizi, M., and Bouhaddi, N., 2007.
\newblock ``{Robustness of mechanical systems against uncertainties}''.
\newblock {\em Finite Elements in Analysis and Design, {\bf 43}}(9), June,
  pp.~715--731.

\bibitem{Ghanmi2011}
Ghanmi, S., Guedri, M., Bouazizi, M.-L., and Bouhaddi, N., 2011.
\newblock ``{Robust multi-objective and multi-level optimization of complex
  mechanical structures}''.
\newblock {\em Mechanical Systems and Signal Processing, {\bf 25}}(7), Oct.,
  pp.~2444--2461.

\bibitem{Forouraghi2000}
Forouraghi, B., 2000.
\newblock ``{A genetic algorithm for multiobjective robust design}''.
\newblock {\em Applied Intelligence, {\bf 12}}(3), pp.~151--161.

\bibitem{Olvander2005}
\"{O}lvander, J., 2005.
\newblock ``{Robustness considerations in multi-objective optimal design}''.
\newblock {\em Journal of Engineering Design, {\bf 16}}(5), Oct., pp.~511--523.

\bibitem{Avigad2005}
Avigad, G., and Moshaiov, A., 2005.
\newblock ``{MOEA-based approach to delayed decisions for robust conceptual
  design}''.
\newblock In Applications of Evolutionary Computing Lecture Notes in Computer
  Science, Springer Berlin / Heidelberg, pp.~584--589.

\bibitem{Xu2011}
Xu, H., Huang, H., and Wang, Z., 2011.
\newblock ``{Research on multi-objective robust design}''.
\newblock In Quality, Reliability, Risk, Maintenance, and Safety Engineering
  (ICQR2MSE), 2011 International Conference on, pp.~885 --890.

\bibitem{Cromvik2011}
Cromvik, C., Lindroth, P., Patriksson, M., and Stromberg., A., 2011.
\newblock {A New Robustness Index for Multi-Objective Optimization based on a
  User Perspective}.
\newblock Tech. rep., Department of Mathematical Sciences, Chalmers University
  of Ttechnology and University of Gothenburg., Gothenburg, Sweden.

\bibitem{Li2008}
Li, M., and Azarm, S., 2008.
\newblock ``{Multiobjective collaborative robust optimization with interval
  uncertainty and interdisciplinary uncertainty propagation}''.
\newblock {\em Journal of Mechanical Design, {\bf 130}}(8), p.~081402.

\bibitem{Du2009}
Du, X., Venigella, P.~K., and Liu, D., 2009.
\newblock ``{Robust mechanism synthesis with random and interval variables}''.
\newblock {\em Mechanism and Machine Theory, {\bf 44}}(7), July,
  pp.~1321--1337.

\bibitem{Fonseca1998}
Fonseca, C., and Fleming, P., 1998.
\newblock ``{Multiobjective optimization and multiple constraint handling with
  evolutionary algorithms. I. A unified formulation}''.
\newblock {\em IEEE Transactions on Systems, Man, and Cybernetics - Part A:
  Systems and Humans, {\bf 28}}(1), pp.~26--37.

\bibitem{Augusto2006}
Augusto, O., Rabeau, S., Depince, P., and Bennis, F., 2006.
\newblock ``{Multi-objective genetic algorithms: A way to improve the
  convergence rate}''.
\newblock {\em Engineering Applications of Artificial Intelligence, {\bf
  19}}(5), Aug., pp.~501--510.

\bibitem{McKay1979}
McKay, M., Beckman, R., and Conover, W., 1979.
\newblock ``{Comparison of three methods for selecting values of input
  variables in the analysis of output from a computer code}''.
\newblock {\em Technometrics, {\bf 21}}, pp.~239--245.

\bibitem{Luo2008a}
Luo, B., 2008.
\newblock ``{Efficient MOEAs with an Adaptive Sampling Technique in Searching
  Robust Optimal Solutions}''.
\newblock In Proceedings of the 7th World Congress on Intelligent Control and
  Automation, no.~863, pp.~117--123.

\bibitem{Petrone2011}
Petrone, G., Nicola, C.~D., Quagliarella, D., Witteveen, J., and Iaccarino, G.,
  2011.
\newblock ``{Wind Turbine Performance Analysis Under Uncertainty}''.
\newblock In 49th AIAA Aerospace Sciences Meeting including the New Horizons
  Forum and Aerospace Exposition, no.~January, pp.~AIAA 2011--544.

\bibitem{Wang2012a}
Wang, W., Caro, S., Bennis, F., and {Salinas Mejia}, O.~R., 2012.
\newblock ``{Optimal design of a simplified morphing blade for fixed-speed
  horizontal axis wind turbines}''.
\newblock In Proceedings of the ASME 2012 International Design Engineering
  Technical Conferences, pp.~DETC2012--70225.

\bibitem{Deb2002}
Deb, K., Pratap, A., and Agarwal, S., 2002.
\newblock ``{A fast and elitist multi-objective genetic algorithm: NSGA-II}''.
\newblock {\em IEEE Transaction on Evolutionary Computation, {\bf 6}}(2),
  pp.~182--197.

\end{thebibliography}
\end{document}